А. Н. ГЕРЕГА

# КОНСТРУКТИВНЫЕ ФРАКТАЛЫ В ТЕОРИИ МНОЖЕСТВ

**А. Н. Герега**

**КОНСТРУКТИВНЫЕ ФРАКТАЛЫ
В ТЕОРИИ МНОЖЕСТВ**

**КАНТОРОВЫ ДИСКОНТИНУУМЫ
И КОНТИНУУМЫ СЕРПИНСКОГО**

*Учебное пособие для студентов физико-математических
и инженерных специальностей высших учебных заведений*

Одесса
«Освита Украины»
2017

Классические геометрические фракталы канторово множество и континуумы Серпинского представлены в пособии как теоретико-множественные объекты. Ряд свойств этих множеств анализируется методами математического и функционального анализа, теории итерированных функций и топологии. В приложении дан обзор некоторых универсальных и специальных размерностей, приведены примеры их использования в физических исследованиях.

Пособие предназначено для студентов старших курсов физико-математических и инженерных факультетов, может быть интересно аспирантам, научным сотрудникам и преподавателям вузов.

# ПРЕДИСЛОВИЕ

Основы фрактальной теории [1-4], опубликованные Бенуа Мандельбротом в 1975-77 годах, за прошедшие сорок лет оформились в мощное научное направление, одним из центров внимания которого стали конструктивные (геометрические) фракталы. В сложившемся в научной и научно-популярной литературе стереотипе их описания основное внимание уделяется изучению техники построения, асимптотическому поведению, расчету размерностей и т.п. Между тем, такие фракталы как канторово множество, ковер и треугольник Серпинского – традиционные объекты исследования теории множеств, топологии, математического и функционального анализа, теории итерированных функций [5-15]. Именно в трактовке математических дисциплин, в которых они были введены в научный обиход и впервые исследовались, эти объекты рассматриваются в настоящем пособии. Такой ракурс возвращает фракталы в круг интересов классической математики, актуализирует возможности аппарата этих дисциплин, провоцирует возникновение удачных аналогий, конструктивных ассоциаций и идей.

В Приложении к пособию рассмотрены некоторые универсальные и специальные размерности, вошедшие в математический аппарат теоретической физики в двадцатом столетии.

# КАНТОРОВЫ ДИСКОНТИНУУМЫ

> *Канторизм – это болезнь, от которой математике еще предстоит избавиться.*
> *Анри Пуанкаре*

**Канторовы множества: определение и свойства**

Генри Д. Смит в 1875 году (или еще раньше) и Георг Кантор в 1883 г. описали ставшее классическим триадное множество исключенных средних, получившее название канторова множества [8].

Множество Кантора может быть построено по несложному алгоритму [8, 13, 16]: отрезок [0,1] делится на три равные части, центральная часть – интервал (1/3, 2/3) – извлекается; причем важно, что именно интервал, т.е. без конечных точек. Затем каждый из оставшихся сегментов делится на три части и удаляются средние. На языке «инициатор-генератор» это значит, что инициатор – это единичный отрезок, а генератор – такой же отрезок, разделенный на 3 равные части, из которого вынута центральная. Процедура повторяется бесконечно. На каждом *n*-м шаге удаляются $N = 2^{n-1}$ интервалов, имеющих длину $r = 3^{-n}$ каждый; таким образом, может быть получена бесконечная последовательность множеств

$C_1$= [0,1/3]∪[2/3,1]; $C_2$= [0,1/9]∪[2/9,1/3]∪[2/3,7/9]∪[8/9,1];

$C_3$= [0,1/27]∪[2/27,1/9]∪[2/9,7/27]∪[8/27,1/3]∪[2/3,19/27]∪
    ∪[20/27,7/9]∪[8/9,25/27]∪[26/27,1]

и т. д., пересечение которых есть триадное канторово множество.

На рис. 1 показаны первые четыре шага построения множества.

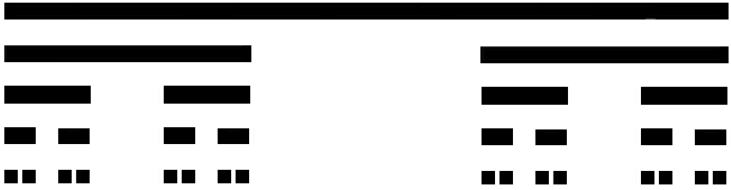

Рис. 1. К построению триадного множества Кантора [16]

Канторово множество точек может быть представлено бесконечными троичными дробями $0,a_1 a_2 a_3 \ldots$, у которых $a_1$, $a_2$, $a_3$,… принимают значения 0 и 2, но не 1. На первом шаге построения множества из отрезка [0, 1] удаляются все точки, у которых первый знак 1, кроме точки $0,1000\ldots = 0,0222\ldots$; на втором шаге из оставшихся удаляются все те, у которых второй знак 1, кроме точек $0,0100\ldots = 0,0022$ и $0,210\ldots = 0,2022$ и т.д. [13]. Понятно, что концы интервалов – троичные дроби $0,a_1 a_2 a_3 \ldots$ все знаки которых, начиная с некоторого, есть только нули или единицы. Для сегментов, получаемых на втором шаге, например, в качестве первого троичного знака нужно взять соответственно 0 или 2, а в качестве остатков – дроби, представляющие соответствующие концы оставшихся сегментов, которые ограничивают интервалы (1/9, 2/9) и (7/9, 8/9), а именно 0,01 и 0,111… и 2,01 и 2,111…[16].

Интересно, что *любое* число из сегмента [0, 2] можно представить в виде суммы двух канторовых чисел [15], значит, и любое действительное число можно «набрать» как сумму некоторого количества канторовых чисел.

Множество не содержит интервалов положительной длины: это видно по тому, что суммарная длина удаленных отрезков равна единице

$$S = 1/3 + 2/3^2 + 2^2/3^3 + \ldots + 2^{n-1}/3^n =$$
$$= (1/3)(1 + 2/3 + (2/3)^2 + \ldots + (2/3)^n) = (1/3)\frac{1}{1-2/3} = 1,$$

или, что эквивалентно, суммарная длина оставшихся на *n*-м сегментов $(2/3)^n$ стремится к нулю при $n \to \infty$. Таким образом, канторово множество содержит бесконечное количество точек и имеет меру Лебега [8, 13]

$$M = 1 - \sum_{n=1}^{\infty} \frac{2^{n-1}}{3^n} = 0.$$

С точки зрения топологии построенное множество – компактное, совершенное и вполне разрывное (дисконтинуум) [8]. Множество называется компактным, если оно ограничено и замкнуто; замкнутым, если содержит все свои граничные элементы; совершенным, если замкнуто и не содержит изолированных точек. Вполне разрывным называется множество, не содержащее ни одного интервала, и, как следствие, не имеющее внутренних точек [8]. (Строгие определения приведены в разделе Справочные материалы).

Наличие у канторова множества таких свойств означает, в частности, и то, что любое компактное, совершенное и вполне разрывное множество можно преобразовать в классическое множество Кантора, причем существует и обратное преобразование. Всякое такое множество принято называть множеством Кантора [8].

Как известно, если исходное множество компактно, обратная функция непрерывна; взаимно однозначная непрерывная функция, обладающая непрерывной обратной, называется гомеоморфизмом или топологическим отображением. В этом случае множество и его прообраз гомеоморфны, или топологически эквивалентны [8]. Сохраняющиеся при гомеоморфизме свойства множеств – компактность, связность, совершенность – называются топологическими инвариантами. Имеет место теорема: свойство быть канторовым множеством является топологическим инвариантом [8].

Как и континуум, множество Кантора несчетно (это можно показать так же, как и несчетность самого континуума) и нигде не плотно, т.е. каждый интервал прямой содержит подинтервал, целиком свободный от точек данного множества (см. Справ. материалы).

Канторово множество самоподобно [15, 17]. Действительно, если в любой строке множества на рис. 1, начиная со второй, отбросить правый сегмент и увеличить левый в три раза, получим предшествующую строку. Таким образом, множество инвариантно по модулю 1 относительно преобразования подобия с коэффициентом 3 [15]. С другой стороны, построенное таким образом множество не вполне самоподобно [17]: для корректной реализации самоподобия его нужно экстраполировать, дополнив исходное множество еще двумя и расположив их в сегментах [1, 2] и [2, 3]. Повторяя эту процедуру неограниченно, получим самоподобное множество на полупрямой [0, ∞].

**Меры и размерности**

Определим меру канторова множества, используя пробную функцию $M_d$, которая в зависимости от выбора её размерности $d$, обращается в нуль или бесконечность при $\delta \to 0$. Введём размерность Хаусдорфа-Безиковича $D_H$, при которой мера $M_d$ изменяет значение с нуля на бесконечность

$$M_d = \sum \gamma(d)\delta^d = \gamma(d)N(\delta)\delta^d \underset{\delta \to 0}{\to} \begin{cases} 0, & \text{при } d > D, \\ \infty, & \text{при } d < D, \end{cases}$$

где $\gamma(d)$ – геометрический коэффициент, зависящий от формы элементов, покрывающих множество [18]. Оценка размерности Хаусдорфа-Безиковича канторова множества [3, 17] дает значение $D_H = \ln 2 / \ln 3$ (см. Приложение).

В [8] доказана теорема: топологическая размерность $D_T$ классического множества Кантора равна нулю. (Отсюда другое название множества – канторова пыль).

Объекты, для которых выполняется неравенство $D_H > D_T$, в работе [3] названы фракталами. Таким образом, выполнение этого неравенства для канторова множества позволяет классифицировать его как фрактальное [3, 17].

Для описания фракталов в той же работе Б. Мандельброт вводит понятие фрактальной размерности множества и полагает ее равной размерности Хаусдорфа-Безиковича [3], следовательно, фрактальная размерность множества Кантора равна

$$D = \ln 2 / \ln 3$$

(см. также Приложение, стр. 64).

Размерность самоподобия $D_S$ канторова множества определяется выражением [17]

$$D_S = -\lim_{l \to 0} \frac{\ln N(l)}{\ln l} = \lim_{n \to \infty} \frac{\ln 2^n}{\ln 3^n} = \frac{\ln 2}{\ln 3} = 0,6309...$$

Понятие размерности самоподобия позволяет целенаправленно строить множества с любой наперед заданной размерностью. Например, для множества Кантора с $D_S = \tfrac{1}{2}$ таковы множества с количеством интервалов $N = 2$, имеющих длину $r = 1/4$, а также $N = 3$ и $r = 1/9$, ... $N = 7$ и $r = 1/49$ и другие, для которых $1/r = N^2$ (рис. 2).

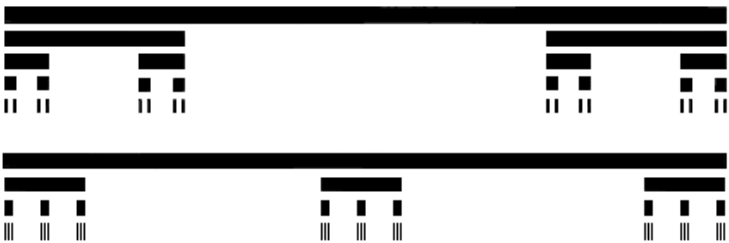

Рис. 2. Два предканторова множества размерности $D_S = \tfrac{1}{2}$: $N = 2$, $r = 1/4$ и $N = 3$, $r = 1/9$ [17]

Рассмотрим вариацию множества Кантора. Пусть исходный единичный отрезок разделен на 5 равных частей. Центральный интервал (2/5, 3/5) – извлекается. Каждый из оставшихся двух сегментов [0, 2/5] и [3/5, 1] делится на 5 равных частей, и в каждом удаляется средний. Процедура повторяется бесконечно, и на каждом $n$-м шаге удаляется $2^{n-1}$ интервалов, имеющих длину $5^{-n}$ каждый; таким образом, может быть получена последовательность множеств, объединение которых – канторово множество с размерностью подобия

$$D_S = \frac{ln\,4}{ln\,5} = 0,86135...$$

Полученное множество не содержит интервалов положительной длины, т.к. суммарная длина удаленных отрезков

$$S = \frac{1}{5} + 2 \cdot \frac{1}{5} \cdot \frac{2}{5} + 4 \cdot \frac{1}{5} \cdot \frac{4}{25} + 8 \cdot \frac{1}{5} \cdot \frac{8}{125} + ... = \frac{1}{5} \cdot \left(1 + \frac{4}{5} + \frac{16}{25} + \frac{64}{125} + ...\right) =$$
$$= \frac{1}{5} \cdot \frac{1}{1 - 4/5} = 1.$$

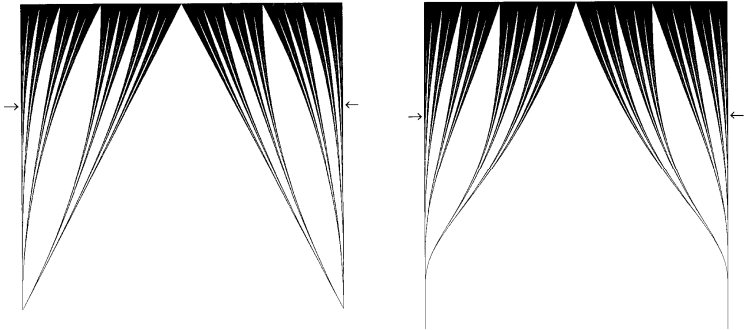

*Любой горизонтальный срез – определенное канторово множество. Горизонтальный отрезок – сегмент* [0, 1]; *по вертикали отложены: слева – значения r, возрастающие сверху вниз от* 0 *до* 1, *справа – размерность D, возрастающая снизу вверх от* 0 *до* 1; *стрелками указаны значения N* = 0,6309 *и r* =1/3 [4].

Рис. 3. Совокупность множеств Кантора в виде «занавеса».

Этот результат легко обобщается для произвольного нечетного числа равных частей исходного отрезка.

Рис. 3 дает представление о форме множеств Кантора с $N = 2$

и различными значениями *r*.

Обобщим алгоритм построения множеств Кантора с заданной размерностью для произвольного соотношения между *N* и *r*.

Пусть, например, нужно построить канторово множество *C*, имеющее размерность самоподобия

$$D_S = \frac{\ln 9}{\ln 10} = 0,9542...,$$

и пусть $C_0$ есть множество всех действительных чисел отрезка [0, 1]. Пусть также множество $C_1$, полученное из $C_0$, содержит числа, в десятичном представлении которых отсутствует одна из цифр, например, тройка [8]. Это значит, что, если разделить $C_1$ на десять равных отрезков, то оно не будет содержать чисел, попадающих в четвертый интервал, граничные точки которого (0,2999…, 0,4000…). Для получения $C_2$ разделим каждый из девяти оставшихся в $C_1$ интервалов на 10 равных частей, и извлечем в них каждый четвертый. Повторяя процедуру бесконечно, получим последовательность множеств $C_1, C_2, C_3, \ldots C_n$, пересечение которых есть искомое множество *C*. Действительно, по построению $C_i$ есть объединение на каждом шаге девяти уменьшенных в десять раз копий самого себя, т.е. *N* = 9, *r* = 1/10 [8].

Используя этот алгоритм, и выбирая для изъятия по одной другие цифры из десятичного представления чисел при новых построениях, можно получить *еще* девять канторовых множеств той же размерности. Понятно, что количество вариантов множеств определяется числом размещений $A_n^m$: изымая по две цифры, получим 90 вариантов множеств Кантора размерности

$$D_S = \frac{\ln 8}{\ln 10} = 0,9031...,$$

по три – 720 множеств размерности

$$D_S = \frac{\ln 7}{\ln 10} = 0,8451...$$

и так далее.

Таким образом, размерность множества $D_S$ и масштаб разбиения $N$ определяют длину интервалов $r = N^{-1/D}$ на каждом шаге, но не их расположение на отрезке [0, 1]. Обсудим возможности количественной оценки различия *двух* множеств одинаковой размерности.

**Упорядоченность и лакунарность**

Рассмотрим меры упорядоченности [19] и лакунарности [3] множеств.

Предложенная в [19] мера упорядоченности базируется на расчете энтропии распределения точек множества. Показано, что мера относительной степени упорядоченности двух равновеликих (для корректности сравнения) множеств есть функционал Ляпунова, известный также как расстояние Кульбака-Лейблера [20, 21] – мера различия вероятностных распределений, и равный

$$\Delta S = -\sum_{i=0}^{255}\left[f_1(i)\ln f_1(i) - f_2(i)\ln f_2(i)\right] = -\sum_{i=0}^{255}\left[f_1(i)\left(\ln f_1(i)/f_2(i)\right)\right],$$

где $f_i$ – функции плотности распределения точек множества. Расчет упорядоченности канторовых множеств на рис. 2 по этому функционалу показывает, что верхнее множество – менее упорядоченное.

Другую возможность количественной оценки дает предложенная Б. Мандельбротом в [3] мера неоднородности множества – лакунарность (от лат. *lacuna* – полость, пустое место, углубление). Он предлагает измерять степень лакунарности канторовой пыли по относительной длине наибольшего удаленного интервала, таким образом, лакунарность триадного канторова множества $L = 1/3$. В двумерных множествах (рис. 4) лакунарность по Мандельброту пропорциональна отношению квадратного корня из площади максимальной полости к ее периметру $\sqrt{S_{max}}/p_{max}$ [3]. Заметим, что в этом соотношении заложена неоднозначность: из элементарной геометрии известно, что площадь изопериметрических фигур зависит от их формы [22] и, следовательно, влияет на величину лакунарности.

Для ковров, изображенных на рис. 4, с учетом, что сторона большей лакуны левого равна 150 единицам длины, а правого – 50, это соотношение дает значения лакунарностей, равные соответственно $L \sim 150/600 = 1/4$ и $L \sim \sqrt{9 \cdot 250}/(9 \cdot 200) = 1/12$. Таким образом, множество, содержащее бóльшие пустоты – более лакунарно [4]; интересно, что при этом степень его упорядоченности меньше, чем у правого ковра. Важно понимать, что лакунарность «не имеет ничего общего с топологией и касается лишь различий фракталов одинаковой размерности» [4].

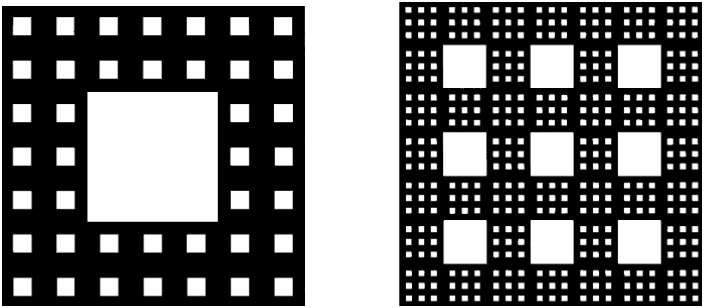

Рис. 4. Два предфрактала различных ковров Серпинского, имеющих одинаковую размерность $D_S = ln\ 40\ /\ ln\ 7 =$
$= 1{,}8957\ldots$ [4].

Вернемся к анализу лакунарности канторовой пыли. Рассмотрим два случая [4]. В первом – остающиеся сегменты канторова множества группируются вблизи границ первоначального отрезка [0, 1], в середине которого возникает большая лакуна размером $1 - N r = 1 - N^{1-1/D}$, что сродни верхнему множеству на рис. 2. В другом случае $N$ сегментов разделены $N - 1$ пустотами, каждая длиной $(1-Nr)(N-1)$, что на рисунке соответствует нижнему множеству: оно содержит небольшие пустоты, и его лакунарность меньше. Таким образом, более упорядоченное множество, как и в случае ковров Серпинского, – менее лакунарно. Последнее утверждение – следствие общего правила.

Для нерегулярных (хаотичных) фракталов необходимы иные способы определения лакунарности. Представляется естественным использовать для такой оценки дисперсию распределения точек множества

$$\langle s^2 \rangle - \langle s \rangle^2 \sim L^{2D},$$

где $s$ – количество элементов в выделяемых для оценки лакунарности кластерах, $D$ – размерность множества. Возможны и другие конструкции.

**Канторова функция**

Обратимся вновь к исходному канторову множеству удаленных средних третей, и построим функцию Кантора – зависимость относительного веса $M(x)$ множества чисел, лежащих на единичном сегменте слева от $x$ [15].

На предфрактале, полученном на первом шаге, по мере того как $x$ проходит сегмент [0, 1/3], $M(x)$ возрастает от 0 до 1/2; на отрезке [1/3, 2/3] – $M(x)$ не изменяется; на [2/3, 1] – функция растет от 1/2 до 1. На следующем шаге построения появятся еще два плато: $M(x) = 1/4$ на [1/9, 2/9] и $M(x) = 3/4$ на [7/9, 8/9]; на третьем шаге: на отрезках [1/27, 2/27], [7/27, 8/27], [19/27, 20/27], [25/27, 26/27] значения $M(x)$ на плато равны соответственно 1/8, 3/8, 5/8, 7/8 (рис. 5). В пределе ступенчатая функция $M(x)$ имеет плато почти всюду, однако, несчетное количество бесконечно малых разрывов позволяет ей изменить свое значение от 0 до 1.

Доопределим функцию по непрерывности: соединим граничные точки соседних плато отрезками прямых, тангенс угла на-клона которых равен $(2r_n)^{-2}$, где $r_n$ – размер сегмента, получаемый на $n$-м шаге построения множества Кантора (рис. 6). Можно показать, что суммарная длина этих отрезков равна единице [14, 23], следовательно, длина «канторовой лестницы» (другое название «чертова лестница») – двум.

Так определенная функция Кантора есть непрерывная, возрастающая, с производной почти всюду равной нулю.

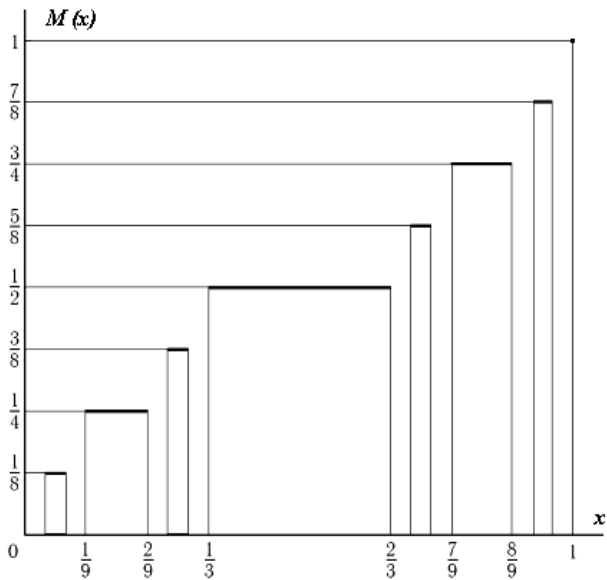

Рис. 5. Функция Кантора на предфрактале второго поколения [15].

Чтобы узнать высоту, на которой расположено плато для любого значения $x$, можно использовать алгоритм, диктуемый аффинной структурой канторовой функции: запишем значение $x$ в виде троичного числа, затем превратим его в двоичную дробь, двигаясь слева направо и заменяя каждую цифру 2 на 1 вплоть до первой единицы [15]. Первую единицу не изменяем, а все стоящие справа от нее цифры заменяем нулями, получая, таким образом, для каждого $x$ единственное значение функции $M(x)$.

Чтобы по значению функции получить значение или интервал значений *x*, запишем $M(x)$ в виде двоичной дроби, заменим все единицы, кроме последней двойками, и каждый нуль, расположенный правее последней единицы, заменим на 1, 2 или 3, получая некоторый интервал [15].

Можно показать [15], что любое значение $M(x)$, представленное в виде дроби с взаимно простыми числителем и знаменателем, которые записаны в десятичной системе счисления, и в которых знаменатель есть *n*-я степень двух, лежат на плато шириной $3^{-n}$. Остальные значения функции в виде бесконечных двоичных дробей соответствуют единственным значениям *x* [15].

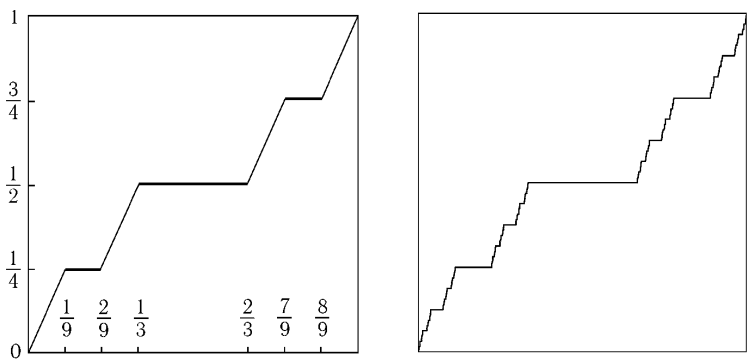

Рис. 6. Второй и более продвинутый этапы построения «чертовой лестницы» [15].

Рассмотрим еще один способ определения канторовой функции [17]. Пусть исходный элемент множества – стержень, сделанный из материала плотности $\rho_0 = 1$. Тогда, при исходной длине стержня, равной единице, его масса $\mu_0 = 1$. Алгоритм построения триадного канторова стержня таков: разрежем стержень на две равные половины $\mu_1 = \mu_2 = 0{,}5$, которые в результате «ковки» укорачиваются до длины $l_1 = l_2 = 1/3$, что приводит к росту плотности до величины $\rho_1 = \mu_1 / l_1 = 3/2$. Бесконечно повторяя эту процедуру на абстрактном стержне, будем получать в $n$-м поколении $N_n = 2^n$ малых стержней длиной $l_n = 3^{-n}$ и массой $\mu_n = 2^{-n}$ [17], при этом суммарная масса стержней остается равной исходной [1]

$$\sum_{n=1}^{N} \mu_n = 1.$$

На рис. 7 показаны первые пять шагов построения триадного канторова стержня, причем высота фрагментов на каждом шаге пропорциональна плотности. Так как масса стержня $n$-го поколения зависит от длины как $\mu_n = l_n^\alpha$, следовательно, его плотность равна

$$\rho_n = \mu_n/l_n = \rho_0 l_n^{\alpha-1}.$$

Здесь $\alpha$ – скейлинговый показатель, известный как показатель Липшица-Гёльдера [15], для такого стержня равный $\alpha = ln2/ln3$.

В данном случае функция Кантора определяет распределение массы канторова стержня [8, 17]. Пусть стержень расположен вдоль оси абсцисс, и начало координат совпадает с левым концом стержня. Тогда масса, содержащаяся в сегменте $[0, x]$ равна

$$M(x) = \int_0^x \rho(t)\,dt.$$

(Функция $\rho(t)$, играющая роль плотности, равна бесконечности во всех точках, образующих канторово множество, и нулю во всех пустых интервалах. Такое поведение плотности позволяет использовать для ее расчета переопределенную $\delta$-функцию Дирака [24]).

Функция $M(x)$ может быть мерой любой величины, определенной на геометрическом множестве.

---

[1] Канторово множество – контрпример к древнегреческой антиномии между дискретным и бесконечно делимым: оно и совершенно разрывно, и бесконечно делимо [15].

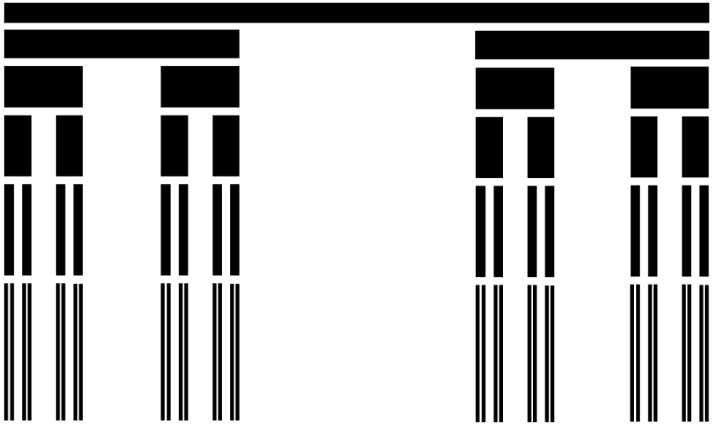

Рис. 7. Пять шагов построения триадного канторова стержня [17]. Высота фрагментов стержня пропорциональна плотности.

**Мультифрактальность**

Еще раз вернемся к исходному канторову множеству удаленных средних, и построим множество Кантора с отличной от нуля мерой Лебега – множество, принадлежащее классу так называемых «толстых» фракталов [15].

Пусть сегмент [0,1] делится на три равные части, и центральная часть удаляется. На второй итерации из обеих оставшихся третей вырежем по 1/9 их длины, т.е. 1/27 единичного отрезка (рис. 8), на третьем – из четырех оставшихся вырежем по 1/81 их длины и т.д. [15]. Таким образом, на каждом шаге удаляются сегменты, относительная длина каждого из которых $(1/3)^{2^k}$, и после $n$ итераций получаем $2^n$ сегментов общей длиной

$$l_n = \prod_{k=0}^{n-1}(1-3^{-2^k}),$$

которая при $n \to \infty$ стремится к $l_\infty = 0{,}585187\ldots$

Рассмотрим два канторовых множества с неравными по длине образующими сегментами, и определим их меру.

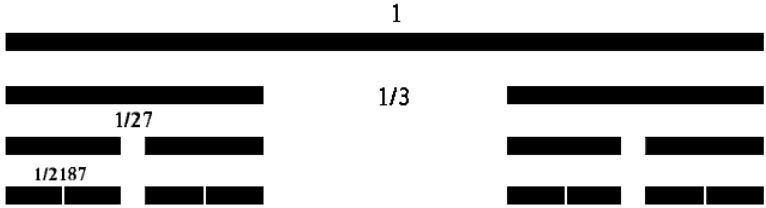

Рис. 8. К построению «толстого» фрактала [15].

Пусть первое множество есть канторов стержень с образующими сегментами относительной длины $l_1 = 1/4$ и $l_2 = 2/5$, которые обладают одинаковыми мерами, равными мере исходного стержня (рис. 9). Определим эту меру.

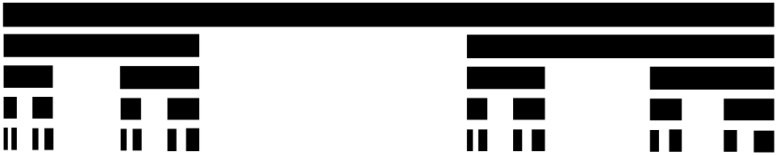

Рис. 9. К построению двухмасштабного канторова множества постоянной меры с $l_1 = 1/4$ и $l_2 = 2/5$ [17].

На $n$-м шаге построения стержень содержит $N = 2^n$ сегментов, причём каждый из $\binom{n}{k} = \dfrac{n!}{k!(n-k)!}$ сегментов имеет соответствующую длину $l_1^k l_2^{n-k}$, $k = 0, 1, \ldots n$. Тогда мера $M_d$ на $n$-м шаге равна

$$M_d = \sum_{i=1}^{N} l_i^d = \sum_{k=0}^{n} \binom{n}{k} l_1^{kd} l_2^{(n-k)d} = \left(l_1^d + l_2^d\right)^n,$$

и остаётся конечной при $n \to \infty$ тогда и только тогда, если $d = D$, где $D$ есть решение уравнения

$$(l_1^D + l_2^D) = 1.$$

Это значит, что так сконструированное множество имеет один скейлинговый показатель – хаусдорфову размерность. Численное решение уравнения для $l_1 = 1/4$ и $l_2 = 2/5$ дает значение размерности $D = 0{,}6110\ldots$ [17].

Необходимость использовать более чем один показатель скейлинга регулярно возникает при исследовании мультифракталов – неоднородных фрактальных объектов, для полного описания которых, в отличие от регулярных фракталов, недостаточно хаусдорфовой размерности, а необходим спектр размерностей, число которых, по существу, бесконечно. Для нас интересно, что мультифрактальные меры могут возникать и при исследовании распределений некой величины на геометрическом (фрактальном) носителе [17].

Рассмотрим канторов стержень с образующими сегментами той же относительной длины $l_1 = 1/4$ и $l_2 = 2/5$, но обладающих мерами $p_1 = 0{,}6$ и $p_2 = 0{,}4$ от меры предшествующего сегмента (рис. 10). Оказывается [4, 16, 17], что в такой постановке множество есть мультифрактал, для адекватного описания которого в данном случае достаточно двух показателей скейлинга: одного для фрактального носителя, другого для вероятностей [15].

Действительно, такой стержень составляют два различных непересекающихся множества, поэтому $d$-меру для определения

$$M_d(q,\delta) = \sum_{i=1}^{N} p_i^q \delta^d = N(q,\delta)\delta^d \to \begin{cases} 0, & \text{при } d > \tau(q) \\ \infty, & \text{при } d < \tau(q) \end{cases}$$

размерности Хаусдорфа-Безиковича канторова стержня можно записать так

$$M_d(q,\delta) = \sum_{k=0}^{n} \binom{n}{k} (p_1^q l_1^d)^k (p_2^q l_2^d)^{n-k} = (p_1^q l_1^d + p_2^q l_2^d)^n.$$

Эта мера остается конечной для $\delta = l_2^n$ при $n \to \infty$ в том и только том случае, когда $d = \tau(q)$, где $\tau(q)$ – решение уравнения

$p_1^q l_1^{\tau(q)} + p_2^q l_2^{\tau(q)} = 1$ (см. Приложение, стр. 68).

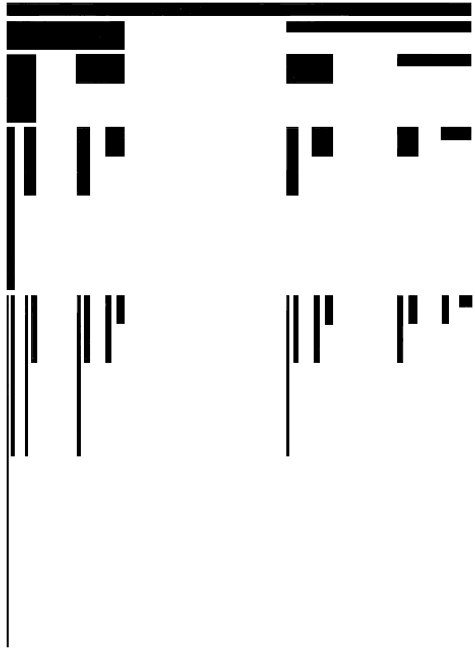

Рис. 10. Двухмасштабное предканторово множество 4-го поколения с сегментами относительной длины $l_1 = 1/4$ и $l_2 = 2/5$ и мерами $p_1 = 0{,}6$ и $p_2 = 0{,}4$ [17].

Уравнение решается численно. При $q = 0$ показатель массы $\tau(q)$ равен $D = 0{,}6110\ldots$; это значение фрактальной размерности. При $q \to \infty$ имеем $p_1^q \gg p_2^q$, следовательно, $\tau(q)$ можно определить из уравнения $p_1^q l_1^{\tau(q)} = 1$ [17]:

$$\tau(q)\big|_{q \to +\infty} = -q\alpha_{min}, \text{ где } \alpha_{min} = \ln p_1 / \ln l_1 = 0{,}3685\ldots.$$

При $q \to -\infty$ верно $p_1^q \ll p_2^q$, тогда

$$\tau(q)\big|_{q \to +\infty} = -q\alpha_{max}, \text{ где } \alpha_{max} = \ln p_2 / \ln l_2 = 1.$$

Авторы [25-27] показали существование связи между показателями массы $\tau(q)$ и размерностями Реньи $D_q$ [28-30] (см. Приложение, с. 70).

Эти размерности определяются соотношением

$$D_q = \tau(q)/(1-q),$$

где

$$\tau(q) = \lim_{\delta \to 0} (\ln Z(q,\delta)/\ln \delta),$$

$$Z(q,\delta) = \sum_{i=1}^{N(\delta)} p_i^q(\delta)$$

– статистическая сумма, а множитель $1/(1-q)$ выбран, чтобы размерность множеств постоянной плотности совпадала с топологической [17].

На рис. 11 показана зависимость фрактальных размерностей спектра Реньи $D_q$ от порядка момента $q$ для рассматриваемого канторова множества [17].

Используя показатели массы $\tau(q)$, можно определить показатель Липшица-Гёльдера $\alpha$

$$\alpha(q) = -\frac{d}{dq}\tau(q)$$

и фрактальные размерности $f(\alpha)$ подмножеств с показателем $\alpha$
$$f(\alpha(q)) = q\,\alpha(q) + \tau(q).$$

Функция $f(\alpha)$ характеризует меры и эквивалентна последовательности показателей $\tau(q)$ (рис. 12) [17].

Для последнего множества введем дополнительное условие. Пусть два интервала, остающиеся от отрезка, образованного на предыдущем шаге итерационной процедуры построения множества, имеют длины, относящиеся как $l_1/l_2$, причём $l_1 + l_2 = 1$, тогда, фигурирующая в выражении для статистической суммы вероятность $p_i^q$ также распадается на два слагаемых так, что
$$Z_{k+1} = (l_1^q + l_2^q) Z_k.$$

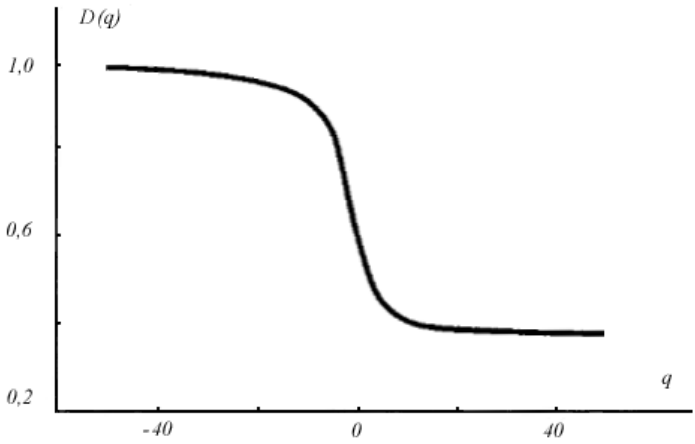

Рис. 11. Спектр фрактальных размерностей как функция порядка момента для триадного канторова множества с $l_1 = 1/4$ и $l_2 = 2/5$, $p_1 = 0{,}6$ и $p_2 = 0{,}4$ [17].

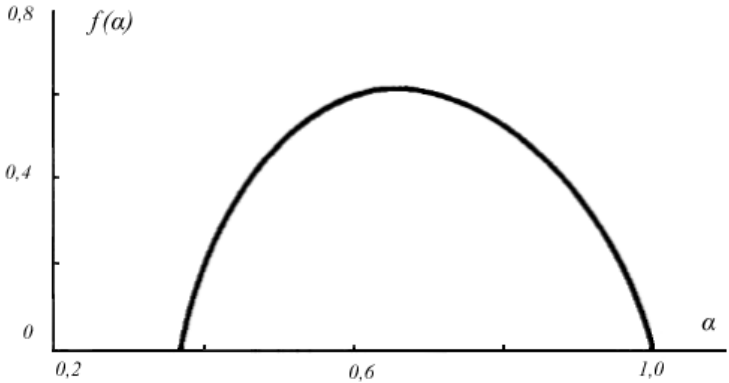

Рис. 12. Зависимость фрактальной размерности $f(\alpha)$ подмножеств канторовой пыли для двухмасштабной фрактальной меры с $l_1 = 1/4$ и $l_2 = 2/5$, $p_1 = 0{,}6$ и $p_2 = 0{,}4$ [17].

С учётом триадности множества, имеем

$$D_q = \frac{1}{1-q} \lim_{k \to \infty} \frac{k\,log(l_1^q + l_2^q)}{-k\,log\,3} = \frac{log_3(l_1^q + l_2^q)}{1-q}.$$

В случае $l_1 = l_2 = 1/2$, все $D_q$ одинаковы и, очевидно, равны хаусдорфовой размерности множества

$$D_q = \frac{log_3(2^{-q} + 2^{-q})}{1-q} = log_3 2,$$

а если $l_1 = 3/4$, $l_2 = 1/4$, канторово множество неоднородно и образует мультифрактал с $D_0 = log_3 2 \approx 0{,}631$, $D_1 \approx 0{,}512$, $D_2 \approx 0{,}428$, $D_\infty \approx 0{,}262$, $D_{-\infty} \approx 1{,}262$ [9].

# КОНТИНУУМЫ СЕРПИНСКОГО

*Фракталы выражаются не в первичных геометрических формах, а в алгоритмах, в наборах математических процедур.*

«Язык фракталов» [31]

*Математическая вселенная населена не только важными видами, но и интересными особями.*

Карл Зигель

**Первый континуум Серпинского**

В 1916 году Вацлав Серпинский построил множество, со временем получившее название первого континуума Серпинского, как пример «кривой, каждая точка которой есть точка ветвления» [6, 32] (см. Справ. материалы). Авторский алгоритм состоит в следующем: в равностороннем треугольнике проводят три средние линии и удаляют внутренние точки образованного ими треугольника. Построение повторяют в каждом из трех оставшихся треугольников, получая 9 треугольников второго поколения, затем 27 – третьего и так далее ad infinitum. Пересечение множеств, полученных на каждом шаге построения, есть треугольник Серпинского (рис. 13).

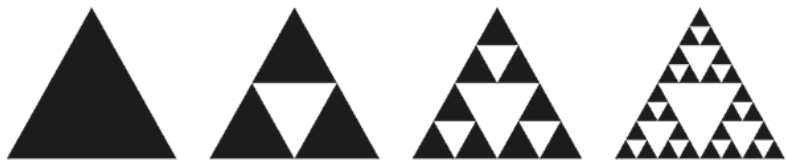

Рис. 13. Три шага построения треугольника Серпинского.

Можно построить треугольник Серпинского по-другому. Как известно, треугольник в евклидовой плоскости – это три точки

(вершины) и три отрезка прямых (стороны) с концами в этих точках [33]. Проведем в треугольнике три средние линии; затем построим средние линии в получившихся малых треугольниках, за исключением перевернутого (рис. 14). Снова, исключая перевернутые, получим 9 треугольников второго поколения. Продолжая построение до бесконечности, получим треугольник Серпинского.

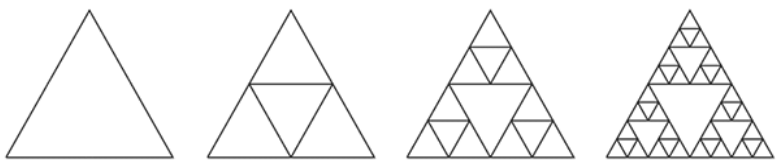

Рис. 14. Три шага другого алгоритма.

Еще один вариант формирования треугольника, предложенный В. Серпинским в той же статье [32], понятен из анализа построения на рис. 15. Выделенные на рисунке ломанные можно положить в основу версии построения, похожей на обычную схему генерации геометрических фракталов с помощью замены частей очередной итерации на масштабированный фрагмент предыдущей. В данном случае – это фрагмент из трех звеньев, который получен в первой итерации, причем откладывать его нужно попеременно вправо и влево. Отметим, что меандры этого построения напоминают кривые типа Пеано-Гильберта.

Таким образом, по построению треугольник Серпинского есть объединение создаваемых на каждом шаге $n \to \infty$ существенно непересекающихся трех копий полученного в предыдущей итерации треугольника, уменьшенных в два раза, следовательно, размерность подобия такого множества равна
$$D = ln3 / ln2 \approx 1,58496250\ldots$$
(В случае пространства произвольной размерности на каждом шаге создается $d+1$ копия, и размерность самоподобия множе-ства равна $D = ln(d+1)/ln\,2$ [3]).

Суммарная площадь удаленных при построении малых треугольников рассчитывается как сумма членов бесконечно убывающей прогрессии

$$S = (1/4)(1 + 3/4 + (3/4)^2 + ... + (3/4)^n) = (1/4)\frac{1}{1 - 3/4} = 1,$$

которая оказывается равной площади исходного.

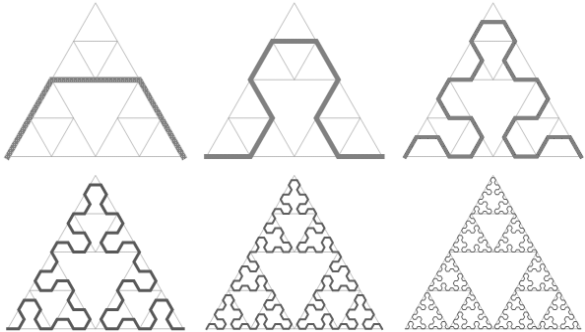

Рис. 15. Вариант построения из статьи В. Серпинского [32].

Представление о нулевой площади треугольника нуждается в уточнении: оно означает, что для любого сколь угодно малого положительного $\varepsilon \to +0$ можно указать фигуру, которая, с одной стороны, содержит треугольник Серпинского, и, с другой, площадь которой не превышает $\varepsilon$ [28].

Определим суммарную величину периметров удаленных в процессе построения треугольников. Пусть сторона исходного треугольника равна единице. Тогда периметр центрального треугольника, полученного в первой итерации, равен 3/2. На второй итерации три малых треугольника имеют суммарный периметр $3 \cdot 3/4 = (3/2)^2$, на $n$-й – их периметр определяется суммой геометрической прогрессии

$$P_n = \sum_{k=1}^{n}\left(\frac{3}{2}\right)^n = \frac{3}{2}\cdot\frac{1-(3/2)^n}{1-(3/2)} = 3\cdot\left(\left(\frac{3}{2}\right)^n - 1\right),$$

и, очевидно, расходится при $n \to \infty$. По аналогии с законом Ричардсона [4] $P_n$ можно представить в виде степенной зависимости

$$P_n = a\cdot\delta^{1-D} = 3\cdot\delta^{1-ln3/ln2},$$

показывающей, что при уменьшении масштаба $\delta$ периметр треугольника неограниченно возрастает.

**Свойства первого континуума**

Пусть $\pi_n$ есть сумма (объединение) всех $3^n$ треугольников ранга $n$. Так как все треугольники замкнуты, то $\pi_n$ – компакт (см. Справ. материалы). Кроме того, это множество, очевидно, связно: любые две его точки могут быть объединены ломаной, лежащей на $\pi_n$. С учетом того, что $\pi_n \supset \pi_{n+1}$, пересечение всех $\pi_n$ есть первый континуум Серпинского [6] (рис. 16).

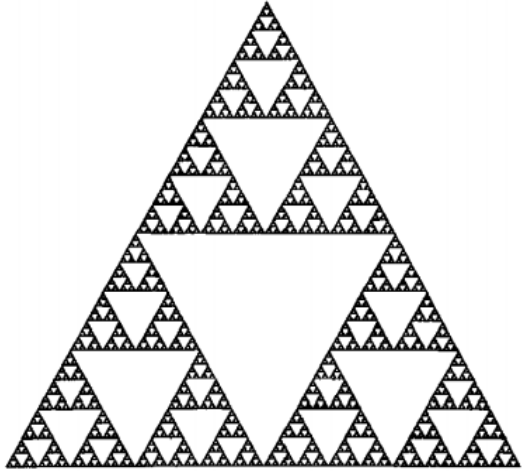

Рис. 16. Континуум «кривая Серпинского» $\pi_6$.

На континууме Серпинского имеется всюду плотная бесконечнозвенная ломаная $L$ – сумма контуров всех треугольников, созданных при построении. Множество $L$ связно как сумма растущей последовательности континуумов $L_n$ – сумм контуров всех треугольников рангов $n$. Интересно, что множество $L$ разбивает свое дополнение на бесконечное количество треугольных областей, обладая при этом нулевой толщиной [8]. Действительно, поскольку по определению Б. Мандельброта [4] для фракталов топологическая размерность $D_T \leq D$, из неравенства $D < 2$ следует, что $D_T$ равна либо нулю, либо единице. Но свойство связности ломаной $L$ оставляет только $D_T = 1$.

Кроме ломаной $L$ на континууме есть еще несчетное множество точек, не лежащих ни на одном из контуров треугольников. Каждая точка этого множества имеет связную окрестность сколь угодно малого диаметра, так что континуум локально связен [6].

**Треугольник Серпинского как аттрактор**

Пусть есть некая динамическая система, эволюция которой описывается набором обыкновенных дифференциальных уравнений

$$\frac{d}{dt}\vec{X}(t) = F(\vec{X}, t),$$

где $\vec{X}(t)$ – вектор в фазовом пространстве $R^n$, $F(\vec{X}, t)$ – векторное поле над этим пространством. Система дифференциальных уравнений такого типа называется потоком в $R^n$ [9, 10].

Пусть каждому состоянию системы соответствует точка в фазовом пространстве решений. Эволюцию системы можно рассматривать как движение изображающей точки по некой траектории в этом пространстве. Известно, что в случае диссипативной динамической системы при $t \to \infty$ все траектории будут сходиться к некоторому компактному множеству в $n$-мерном фазовом пространстве, имеющему нулевой объем, инвариантному относительно действия потока и называемому аттрактором сис-темы [9, 10].

Проведем аналогию между поведением диссипативных систем и алгоритмами построения первого континуума Серпинского.

Рассмотрим систему итерируемых сжимающих функций, реализующих аффинное преобразование точечных множеств, аттрактор которых – треугольник Серпинского

$$f_1\begin{pmatrix}x_1\\x_2\end{pmatrix}=\begin{pmatrix}1/2 & 0\\0 & 1/2\end{pmatrix}\begin{pmatrix}x_1\\x_2\end{pmatrix}$$

$$f_2\begin{pmatrix}x_1\\x_2\end{pmatrix}=\begin{pmatrix}1/2 & 0\\0 & 1/2\end{pmatrix}\begin{pmatrix}x_1\\x_2\end{pmatrix}+\begin{pmatrix}1/2\\0\end{pmatrix} \quad (\star)$$

$$f_3\begin{pmatrix}x_1\\x_2\end{pmatrix}=\begin{pmatrix}1/2 & 0\\0 & 1/2\end{pmatrix}\begin{pmatrix}x_1\\x_2\end{pmatrix}+\begin{pmatrix}1/4\\\sqrt{3}/4\end{pmatrix}.$$

Применим эту систему уравнений к первому треугольнику на рис. 14, предполагая, что его сторона равна единице. Результатом действия каждой из функций $f_1 \div f_3$ будут соответственно треугольники, изображенные на рис. 17 [8, 29].

Произвольное количество применений функций $f_1 \div f_3$ в различной последовательности приведет к образованию всей совокупности треугольников различных поколений континуума Серпинского. На рис. 18 показан треугольник Серпинского, возникающий после четырех итераций [8, 29].

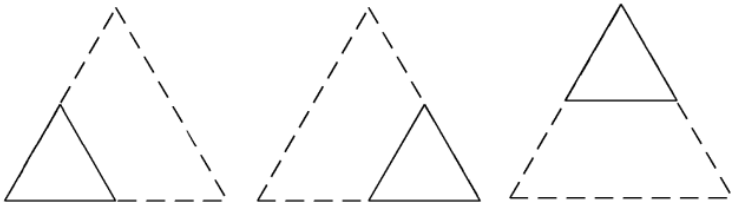

Рис. 17. Результаты действия функций системы ($\star$) на исходный треугольник.

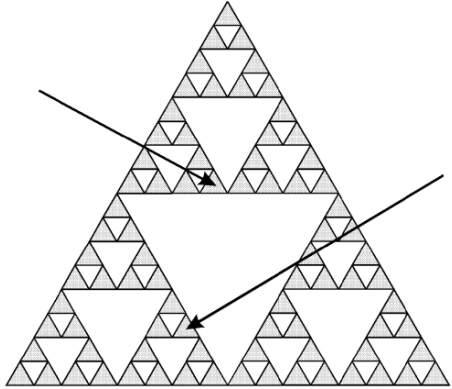

Рис. 18. Треугольник, образованный действием последовательности $f_3\,f_1\,f_2\,f_2$, указан верхней стрелкой, $f_1\,f_2\,f_3\,f_2$ – нижней.

Интересно, что для получения такого же предельного результата можно воздействовать функциями (⋆) на исходную геометрическую фигуру произвольной формы. Например, это может быть квадрат (рис. 19): на каждом шаге уменьшаясь в размерах в два раза и утраиваясь в количестве, эти фигуры образуют треугольник Серпинского – аттрактор системы (⋆) трех линейных преобразований $f_1 \div f_3$.

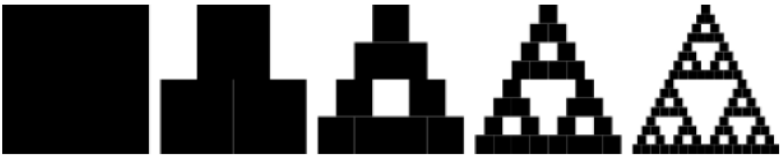

Рис. 19. Вариант построения множества Серпинского с использованием системы итерируемых функций (⋆) [8, 29].

Поскольку треугольник Серпинского – аттрактор, то процесс его построения можно начать даже с единственной точки [29].

Действительно, пусть есть произвольная точка в равностороннем треугольнике. Случайным образом выберем одну из вершин треугольника, и соединим их отрезком прямой. Середина отрезка есть исходная точка для очередной итерации (рис. 20). Этот алгоритм, предложенный Майклом Барнсли [34], по сути, описывает поток, сходящийся к аттрактору в виде треугольного континуума Серпинского. Оказывается, что и при случайном выборе последовательности преобразований (⋆), описывающих движение изображающей точки в фазовом пространстве, при количестве итераций, стремящихся к бесконечности, ее траектория сколь угодно близко подойдет к каждой точке этого фрактального множества. Действительно, пусть, например, начальная точка расположена в центре наибольшего исключенного треугольника (рис. 21); на следующем шаге она – в центре одного из трех треугольников поменьше, т.к. эти треугольники представляют собой геометрическое место точек, которые находятся на половине расстояния до соответствующих вершин от точек большого центрального треугольника. При следующей итерации точка попадает в центр еще меньшего исключенного треугольника и т.д. [29].

Еще один вариант алгоритма построения треугольника – так называемая «игра сэра Пинского» [15]: произвольно выбранную в треугольнике точку соединяют с его ближайшей вершиной; затем от этой вершины откладывают отрезок прямой, проходящий через начальную точку, и равный удвоенному расстоянию от вершины до нее; конец этого отрезка определяет положение следующей точки, и так далее (рис. 22).

Треугольник содержит два типа геометрических мест точек: принадлежащие одному из них точки, будучи выбранными в качестве начальных в «игре», приводят к построению треугольника Серпинского, относящиеся к другому – выводят последовательность за пределы треугольника. На рис. 10 видно, что положение точки 3 спровоцировала точка, находящаяся в белом треугольнике. Его прообразы – три перевернутых треугольника, вдвое меньшие и расположенные по одному в каждом из трех темных. Прообразами этих трех есть девять перевернутых треугольников с вдвое меньшими сторонами, вырезанными из середины девяти оставшихся темных треугольников, линейные размеры каждого из которых в четыре раза меньше размеров исходного [15]. Таким образом, оба набора точек образуют инвариантные множества отображений.

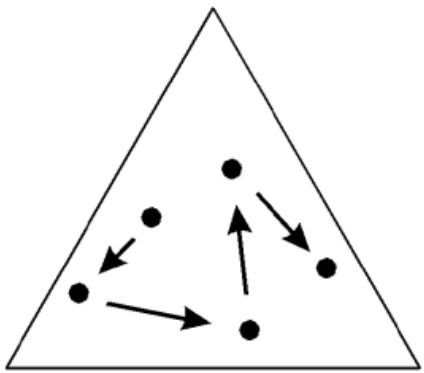

Рис. 20. Иллюстрация к построению треугольника Серпинского по алгоритму М. Барнсли [29].

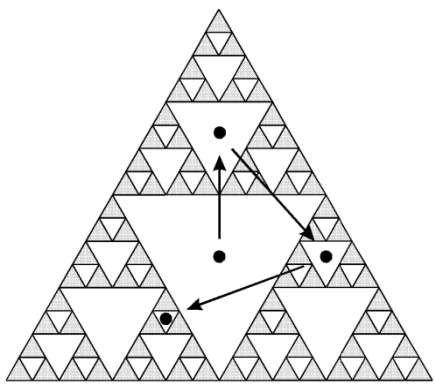

Рис. 21. Треугольник Серпинского как аттрактор случайной последовательности итераций системы функций (⋆).

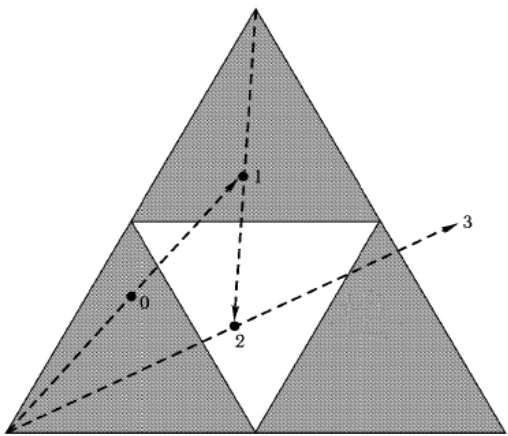

Рис. 22. К алгоритму «сэра Пинского».

Определим точки обоих типов арифметически. Введем систему координат, имеющую симметрию треугольника: оси направим по биссектрисам углов вовне треугольника, вершинам поставим в соответствие единичные значения соответствующих координат, противоположным сторонам – нулевые значения тех же координат [15]. Понятно, что для плоскости набор из трех координат избыточен, но их нельзя выбирать независимо, поскольку они связаны соотношением $x + y + z = 1$, $x > 0$, $y > 0$, $z > 0$.

Можно показать, что для точек, принадлежащих треугольнику Серпинского, на каждом месте двоичной записи координат $x$, $y$, $z$ должны стоять одна единица и два нуля [15]. Три нуля на первом месте после двоичной запятой означают, что значения ни одной из координат не превосходят 1/2; геометрически это значит, что соответствующая точка находится в центральном перевернутом треугольнике, который не принадлежит множеству Серпинского. В общем случае три нуля на $n$-м месте после двоичной запятой соответствует извлечению $3^{n-1}$ перевернутых треугольников со сторонами длиной $2^{-n}$ из $3^{n-1}$ неперевернутых треугольников, оставшихся после $n-1$ операций вырезания центральных перевёрнутых треугольников.

**Пирамида Серпинского**

На рис. 23 – трехмерное обобщение треугольного множества Серпинского. Построение такой пирамиды начинается с правильного тетраэдра, из которого вырезают октаэдр, оставляя в нем четыре малых тетраэдра, с вдвое меньшими гранями (см. рис. 24). Алгоритм повторяют с каждым из малых тетраэдров ad infinitum.

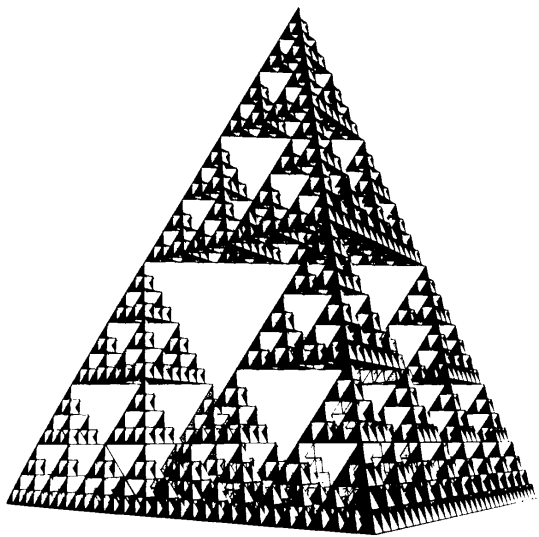

Рис. 23. Пирамида Серпинского [4].

Размерность пирамиды, таким образом, определяется количеством $N = 4$ малых тетраэдров в два раза меньших линейных размеров $r = 1/2$, порождаемых каждым из тетраэдров предшествующей итерации

$$D = \frac{ln4}{ln2} = 2.$$

Это равенство показывает, что алгоритм построения пирамиды Серпинского приводит при $n \to \infty$ к полному изъятию объема исходного тетраэдра.

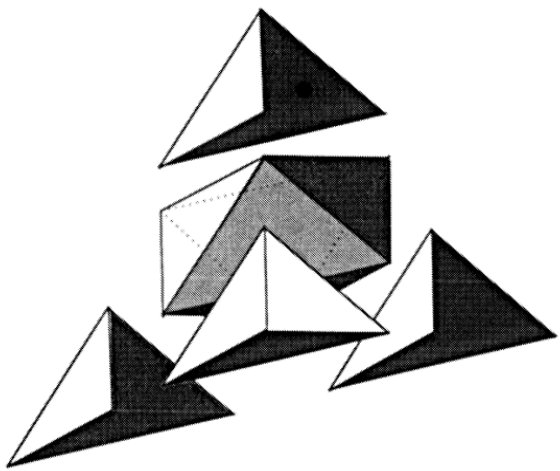

Рис. 24. Центральный октаэдр и малые тетраэдры первого шага построения пирамиды Серпинского [35].

**Второй континуум Серпинского: построение и свойства**

Рассмотрим классическое множество Кантора, обобщенное на случай бóльшего числа измерений.

В двумерном случае возможны два обобщения. Первое – декартово произведение канторова множества на себя $C \times C$. Полученное двумерное множество может быть покрыто $N(r) = 4^n$ квадратами со сторонами длиной $r = 1/3^n$ (рис. 25) и, следовательно, имеет размерность вдвое большую размерности одномерного

$$D(C \times C) = \lim_{n \to \infty} \frac{ln\, 4^n}{ln\, 3^n} = 2 \cdot \frac{ln\, 2}{ln\, 3} \approx 1,2618595...$$

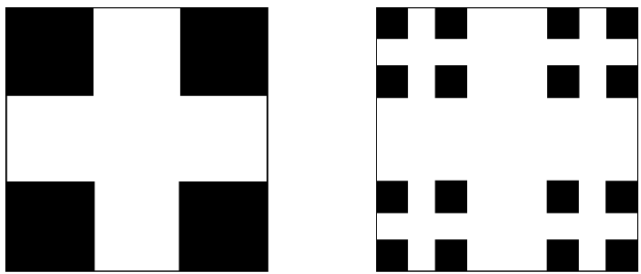

Рис. 25. Предфракталы двумерной канторовой пыли: первый и второй шаг построения.

Второе обобщение – это дополнение к декартову произведению дополнений канторова множества на себя, т.е. $(C' \times C')'$. Это множество допускает удобное рекурсивное построение: инициатор – единичный квадрат, генератор – единичный квадрат, разделенный на 9 малых квадратов со стороной 1/3 исходного, из которого вынут центральный квадрат. На втором шаге процедуру повторяют: в каждом из 8 оставшихся малых квадратов извлекаются центральные со стороной 1/9 исходного и т.д. (рис. 26). Пересечение полученных таким образом после бесконечного количества итераций множеств есть второй континуум Серпинского [12, 36]. Это множество называют также канторов ковер или ковер Серпинского.

Так как множество самоподобно, для определения его размерности достаточно определиться с генератором: для полного покрытия ковра достаточно восемь квадратов со стороной 1/3. Тогда размерность равна

$$D((C' \times C')') = \ln 8 / \ln 3 \approx 1{,}892789\ldots$$

В общем случае хаусдорфова размерность декартова произведения равна сумме размерностей множеств-сомножителей, но проецирование вложенных в многомерные евклидовы пространства фрактальных множеств на евклидовы пространства меньших размерностей дает множества, хаусдорфовы размерности которых зависят от направления проецирования [15]. Например, пусть есть одномерное множество канторовского типа $C_4$, получаемое при удалении из единичного сегмента центральных четвертей, тогда множество ($C_4 \times C_4 \times C_4$) имеет размерность Хаусдорфа

$$D = 3 \cdot \ln 2 / \ln(8/3) = 2{,}120085\ldots$$

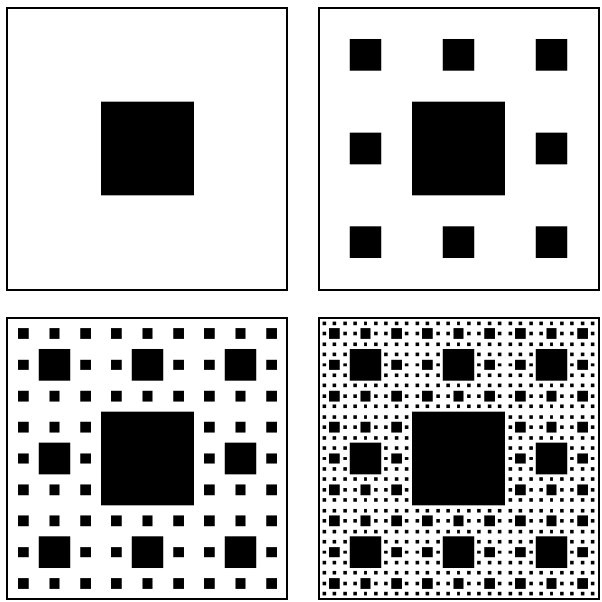

Рис. 26. Предфракталы ковра Серпинского: первые четыре шага.

Проецирование ($C_4 \times C_4 \times C_4$) вдоль одной из трех координатных осей дает множество ($C_4 \times C_4$) с размерностью

$$D = 2 \cdot \ln 2 / \ln(8/3) = 1{,}41339018\ldots$$

Проекции этого же множества на другие направления могут быть двумерными множествами, содержащими связные элементы [15].

Двумерные множества Кантора, как и классические триадные, можно строить с наперед заданной размерностью [8]. Например, несложно получить двумерное множество размерности 1: пусть исходный квадрат на каждом шаге заменяется четырьмя меньшими с $r = 1/4$ (рис. 27), тогда, очевидно, размерность $D = log\,4 / log\,4 = 1$.

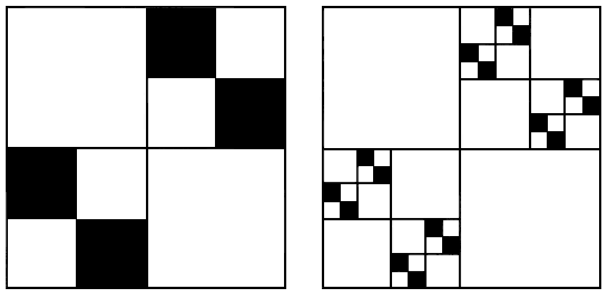

Рис. 27. К построению двумерного множества Кантора единичной размерности: предмножества первого и второго поколения.

Нетроичные квадратные центрированные фрактальные ковры можно построить, положив $r = 1/b$, где $b$ – целое число, большее 3. В качестве тремы (вырезаемой части) возьмем квадрат со стороной $1 - 2r$ с центром в середине исходного квадрата, а в качестве генератора – узкое кольцо из $4 \cdot (b - 1)$ квадратов со стороной $r$. Размерность такого ковра определяется выражением

$$D = ln\left[4(b-1)\right]/ln\,b.$$

Если взять нечетное целое $b > 3$, в качестве тремы – один подквадрат со стороной $r$ и с центром в той же точке, что и центр исходного квадрата, а в качестве генератора – широкое кольцо из $(b^3 - 1)$ малых квадратов, то получится ковер с размерностью

$$D = ln\left(b^3 - 1\right)/ln\,b.$$

Таким образом, в центрированных коврах можно получить сколь угодно близкое приближение к любому значению $D$ в интервале от 1 до 2.

Нецентрированные ковры получаются при $b \leq 2$. Например, при $b = 2$ и $N = 3$ можно разместить трему, состоящую из одного подквадрата, в правом верхнем подквадрате. Соответствующее предельное множество оказывается треугольником Серпинского, построенным из треугольника, образующего левую нижнюю половину квадрата.

**Индекс ветвления ковра Серпинского**

В топологии пользуются понятием линии, введенным в 1921 году П.С. Урысоном, и «являющимся наиболее общим (но не чрезмерно)» [12].

По Урысону линией называется одномерный континуум, т.е. связное компактное метрическое пространство, каждая точка которого обладает сколь угодно малой окрестностью с границей размерности нуль [12]. Другими словами, при любом $\varepsilon > 0$ пространство может быть представлено в виде суммы конечного числа замкнутых множеств диаметра, меньшего $\varepsilon$, обладающих тем свойством, что никакие три из этих множеств не имеют общей точки.

Ковер Серпинского удовлетворяет урысоновскому определению линии, так что всякая канторова кривая является и линией в смысле Урысона [37]. И обратно, если плоский континуум является линией в смысле Урысона, то он есть канторова кривая. Ковер Серпинского – локально связный континуум и потому может быть получен как непрерывный образ отрезка. Любая канторова кривая может быть топологически вложена в ковер Серпинского, т.е. в нем содержится континуум, гомеоморфный линии [12].

Известно, что единственными собственно континуумами (т.е. континуумами, содержащими более одной точки), лежащими на прямой, являются сегменты [6]. Так как при непрерывных отображениях сохраняется и связность, и компактность, то непрерывный образ всякого континуума есть континуум. Отсюда, в частности, следует, что непрерывный образ прямолинейного сегмента есть континуум. Поэтому всякая система $n$ непрерывных функций

$$x_i = x_i(t), i = 1, 2, \ldots, n,$$

заданных на сегменте $0 \leq t \leq 1$, определяет в $n$-мерном пространстве некоторый континуум, являющий в силу самих уравнений непрерывным образом отрезка [0, 1], который называется обычной непрерывной кривой в $n$-мерном пространстве, а система уравнений – параметрическим представлением этой кривой [6]. Такое определение непрерывной кривой даже в двумерном пространстве может включать в себя геометрические образы, вовсе непохожие на то, что мы привыкли называть «линиями»: например, треугольник и квадрат в смысле приведенного определения есть непрерывные кривые. Поэтому континуумы, являющиеся непрерывными образами прямолинейного сегмента, принято называть не кривыми, а жордановыми континуумами [6].

Ковер Серпинского нигде не плотен на плоскости и содержит топологический образ всякого континуума, лежащего в плоскости и нигде не плотного в ней. Такие континуумы называются (плоскими) канторовыми кривыми. Причем, канторова кривая может не быть жордановой и, обратно, жорданов континуум, например квадрат, может не быть канторовой кривой [6].

В исследовании линий важную роль играет понятия индекса ветвления [12] (см. Справочные материалы). Точки линии можно классифицировать как обладающие натуральным, неограниченным, счетным или континуальным индексом ветвления. В ковре Серпинского малые квадраты (ячейки), полученные на произвольном шаге итерационной процедуры построения, счи-таются соединенными, если соприкасаются сторонами; другими словами, ковер Серпинского обладает бесконечной разветвленностью, и задача разделения его на части может быть решена удалением бесконечного (счетного) множества точек. Таким образом, с топологической точки зрения ковер Серпинского – линия с континуальным индексом ветвления во всех своих точках.

Некоторые свойства фракталов с конечной и бесконечной разветвленностью существенно различны [38]. «Наиболее интересным для нас свойством таких фрактальных решеток является то, что в отличие от решеток с конечной разветвленностью, на которых путь протекания разрушается при выбрасывании конечного числа узлов, на этих существует самый настоящий перколяционный переход» [38] (см. Справ. материалы, стр. 51). Параметры перколяционного перехода на ковре Серпинского изучены в [39].

Рассмотрим модификацию ковра Серпинского: пусть соединенными считаются клетки либо соприкасающиеся сторонами, либо имеющие общую вершину. Будем называть такой аналог известного фрактала – ковром Серпинского с гибридной разветвленностью [40]. Понятно, что модификация правил образования связности приводит к изменению перколяционных параметров бесконечного кластера ячеек ковра.

По описанному алгоритму разделим любую ячейку ковра Серпинского произвольного шага итерации на 9 клеток и удалим среднюю. Определим вероятность $p'$ принадлежности ячейки перколяционному кластеру на ковре, т.е. вероятность того, что через ячейку можно «протечь» по составляющим ее клеткам, каждая из которых входит в бесконечный кластер с вероятностью $p$. Так как ренормгрупповое преобразование [41] должно в нашем случае отражать факт наличия связности, количество подходящих комбинаций в расположении клеток в ячейке будет меньше комбинаторного. С учетом этого реном-преобразование для ковра с гибридной разветвленностью имеет вид

$$p' = R(p) = p^8 + 8p^7(1-p) + 27p^6(1-p)^2 +$$
$$+ 44p^5(1-p)^3 + 38p^4(1-p)^4 + 8p^3(1-p)^5,$$

с нетривиальной неподвижной точкой $p_c = 0{,}5093$, определяющей порог протекания.

Индекс длины корреляции перколяционной системы, может быть найден из соотношения $\nu = \ln b / \ln \lambda = 1.801$, где $b = 3$ – количество клеток вдоль стороны ячейки, $\lambda = (dR/dp) \mid p=p_c$. Критический показатель параметра порядка $\beta$ определяется из равенства $D = d - \beta/\nu$, где аппроксимацией размерности $D$ перколяционного кластера служит размерность ковра Серпинского; при размерности пространства $d = 2$ величина $\beta = 0{,}193$. (Для верификации полученных значений: в случае стандартного ковра Серпинского по нашим данным $\nu = 2{,}194$ и $\beta = 0{,}234$, а по результатам работы [38, 39] $\nu = 2{,}13$ и $\beta = 0{,}27$).

Другие критические показатели могут быть определены из системы равенств двухпоказательного скейлинга [38]: индекс средней длины конечного кластера $\gamma = \nu d - 2\beta = 3{,}216$; критический показатель аналога теплоемкости $\alpha = 2 - \nu d = -1{,}602$; определяющий наибольший размер конечных кластеров индекс $\Delta = \nu d - \beta = 1{,}809$.

**Трехмерные аналоги**

В трехмерном пространстве существует два обобщения канторова множества: «канторов сыр» и губка Менгера.

Множество $(C' \times C' \times C')'$ называется «канторовым сыром»; его генератор может быть заполнен 26 кубиками с длиной ребра 1/3. Тогда размерность

$$D((C' \times C' \times C')') = \ln 26 / \ln 3 \approx 2{,}965647\ldots,$$

что очень близко к трем, т.е. «представляет собой достаточно сплошную конструкцию с изолированными пустотами» [15].

Обобщение канторова сыра на $k$ измерений имеет размерность

$$D((C'\times C'...\times C')') = \frac{ln(3^k-1)}{ln3} \approx \frac{1}{ln3}\left(ln3^k + 1/3^k\right) \approx k - 1/(3^k ln3),$$

что, как правило, несколько меньше размерности вложения $k$ [15].

Всякая линия в смысле Урысона гомеоморфна некоторому подмножеству трехмерного евклидова пространства (теорема Менгера) [42]. Континуум $M$, обладающий свойством, что какова бы ни была линия $C$, в нем найдется подконтинуум $C'$, гомеоморфный континууму $C$, строится следующим образом. Куб $K$ с ребром 1 делится плоскостями, параллельными его граням, на 27 равных кубов. Из куба $K$ удаляется центральный и все прилежащие по двумерным граням кубы первого поколения. Получается множество $K_1$, состоящее из 20 замкнутых кубов. Повторяя процедуру с каждым из кубов $K_1$, получим множество $K_2$, состоящее из 400 кубов второго поколения (рис. 28). Продолжая процесс бесконечно, получим последовательность континуумов $K_1 \supset K_2 \supset K_3 \supset ...$, пересечение которых есть одномерный континуум $M$, называемый универсальной кривой Менгера [12, 42].

Пустоты в ней представляют собой открытые каналы, насквозь пронизывающие исходный куб, и каждая грань исходного куба выглядит как ковёр Серпинского. У губки нулевой объем (так как на каждом шаге он умножается на 20/27), бесконечно большая площадь; топологически губка есть универсальная кривая: любая кривая в трехмерном пространстве гомеоморфна некоторому подмножеству губки Менгера. Размерность губки равна

$$D = ln20/ln3 \approx 2,726833...$$

Дополнение губки имеет конечную меру Лебега и размерность, равную трем (рис. 28).

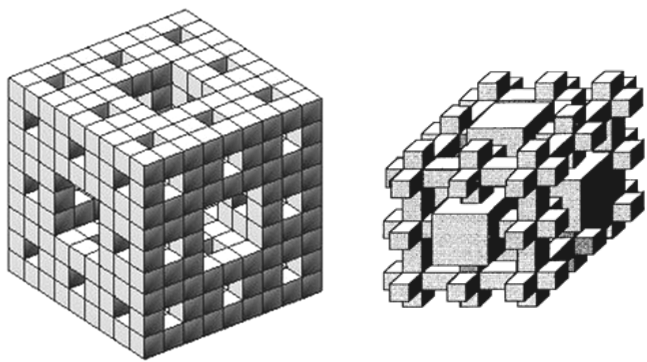

Рис. 28. Универсальная кривая (губка) Менгера и ее дополнение на втором шаге построения.

Теоретико-множественное описание позволяет показать симметрию губки Менгера [15]. Пусть $X$, $Y$, $Z$ – канторовы множества по соответствующим осям. Выражение $X' \cap Y'$ описывает центральный вырезанный квадрат в плоскости $xy$, а запись

$$(X' \cap Y') \cup (X' \cap Z') \cup (Z' \cap Y')$$

– множество всех пустот генератора множества. Дополнение этого множества

$$(X \cup Y) \cap (X \cup Z) \cap (Z \cup Y)$$

есть генератор губки Менгера [15, 43].

Двумерных аналогов губка Менгера не имеет.

## СПРАВОЧНЫЕ МАТЕРИАЛЫ
### О классификации множеств

Пусть $\boldsymbol{R}^n$ – $n$-мерное векторное пространство. Диаметром множества $A \subset \boldsymbol{R}^n$ называется величина

$$\delta(A) = sup\,\{\|\boldsymbol{x} - \boldsymbol{y}\|_2 : \boldsymbol{x}, \boldsymbol{y} \in A\},$$

где $\|\boldsymbol{x} - \boldsymbol{y}\|_2$ – евклидово расстояние между векторами $\boldsymbol{x}$ и $\boldsymbol{y}$; индекс 2 означает, что расстояние равно корню 2-ой степени из суммы абсолютных разностей пар значений, взятых во второй степени. Множество $A \subset \boldsymbol{R}^n$ называется **ограниченным**, если оно имеет конечный диаметр, то есть $\delta(A) < \infty$ [8].

Последовательность $\{\boldsymbol{x}_n\}_{n=1}^{\infty}$ называется сходящейся, то есть имеющей предел $\lim\limits_{n\to\infty}\boldsymbol{x_n} = \boldsymbol{x}$, если для каждого $\varepsilon > 0$ существует такой номер $N$, что при $n > N$ выполняется неравенство $\|\boldsymbol{x} - \boldsymbol{y}\|_2 < \varepsilon$, или, по-другому, $\lim\limits_{n\to\infty}\|\boldsymbol{x_n} - \boldsymbol{x}\|_2 = 0$.

Множество $A \subset \boldsymbol{R}^n$ называется **замкнутым**, если для любой последовательности $\{\boldsymbol{x}_n\}_{n=1}^{\infty}$ точек из $A$, сходящихся к $\boldsymbol{x}$, ее предел также принадлежит множеству $A$: $\boldsymbol{x} \in A$.

Подмножество пространства $\boldsymbol{R}^n$ с евклидовой метрикой **компактно** в том и только в том случае, если оно замкнуто и ограничено [8]. Такое определение компактности множества можно использовать только для подмножеств пространства $\boldsymbol{R}^n$. В случае произвольного метрического пространства множество называется компактным, если из каждой последовательности точек $\{\boldsymbol{x}_n\}_{n=1}^{\infty}$ из множества $X$ можно выделить подпоследовательность $\{\boldsymbol{x}_{n_k}\}_{k=1}^{\infty}$, сходящуюся к некоторой точке $x \in X$.

Эти определения компактности эквивалентны, если $X \subset \boldsymbol{R}^n$.

Точка $x$ множества $A$ есть изолированная точка этого множества, если у нее есть окрестность, не содержащая точек множества $A$. Множество называется **совершенным**, если оно замкнуто и не содержит изолированных точек [8].

Множество $A$ есть **связное** множество, если его нельзя представить в виде объединения двух непустых множеств. Говорят, что множество $A$ **вполне разрывно**, или **вполне несвязно**, если наибольшие связные подмножества $A$ представляют собой одноточечные множества, т.е., если все компоненты $A$ – одинокие точки [8].

Множество $E$ называется плотным на $M$, если каждая точка множества $M$ является предельной точкой множества $E$, т.е. в любой окрестности имеются точки, принадлежащие $E$. Плотные множества на всей прямой называются всюду плотными. Множество называется нигде не плотным (на прямой), если оно неплотно ни на каком интервале, иными словами, если каждый интервал прямой содержит подинтервал, целиком свободный от точек данного множества. Аналогично определяются множества, нигде не плотные на плоскости или, вообще, в произвольном топологическом пространстве. Для того чтобы замкнутое множество было нигде не плотным, необходимо и достаточно, чтобы его дополнение было всюду плотно [5].

**О несчетности континуума**

Совокупность всех точек интервала будем называть континуумом. Покажем, что множество всех действительных чисел, заключенных между 0 и 1, т.е. континуум, несчетен [13].

Пусть все числа, заключенные в интервале (0,1), можно расположить в последовательности $x_1$, $x_2$, $x_3$ …, и пусть каждое из них представлено в виде бесконечной десятичной дроби. «Конечные дроби» дополним до бесконечных нулями, избавляясь от дробей, оканчивающихся последовательностью девяток. Построим десятичную дробь λ такую, что ее *n*-й десятичный знак на единицу больше *n*-го десятичного знака дроби $x_m$, если последний равен 0, 1, 2, …, 7; если же он равен 8 или 9, то *n*-й знак дроби λ равен нулю. Это описание определяет λ полностью, причем она не оканчивается последовательностью девяток. Но таким образом λ есть число из интервала (0,1), отличное от всех чисел $x_n$, и это противоречит исходному предположению, что последовательность $x_1, x_2, x_3$ … содержит все действительные числа, заключенные между 0 и 1.

**Условие Липшица-Гёльдера**

Это условие – неравенство, в котором приращение функции оценивается по приращению аргумента [12]. Функция *f* (*x*), определенная в некоторой области *n*-мерного евклидова пространства $R^n$, удовлетворяет в точке $y \in R^n$ условию Гёльдера с показателем *α* (0 < *α* ≤ 1) и коэффициентом *A* (*y*), если

$$|f(x)-f(y)| \leq A(y)|x-y|^\alpha$$

для всех $x \in R^n$, достаточно близких к *y*. Условие, по сути, означает, что функция изменяется не быстрее, чем некоторая прямая с угловым коэффициентом *A*(*y*). Функция, удовлетворяющая условию Липшица-Гёльдера – ограничена, равномерно непрерывна и дифференцируема.

**Индекс ветвления**

Топологическое понятие индекс (степень) ветвления независимо определили в начале 20-х годов двадцатого века Павел Урысон и Карл Менгер [42, 44]. Индекс ветвления определяется наименьшим количеством точек сечения, которое позволяет разъединить множество.

Линия в точке *x* имеет индекс ветвления *m*, если каково бы ни было число ε > 0, существует открытое множество диаметра, меньшего, чем ε, содержащее точку *x*, граница которого есть множество мощности, не превосходящей *m*, но для достаточного малого ε′ > 0 граница всякого открытого множества, содержащего точку *x*, диаметр которого меньше ε′, имеет мощность, не меньшую, чем *m* [12].

Точка *x* линии имеет неограниченный индекс ветвления, если каково бы ни было число ε > 0, существует открытое множество, содержащее эту точку, с диаметром, меньшим, чем ε, граница которого состоит из конечного множества точек; но каково бы ни было натуральное *n*, найдется такое $ε_n$ > 0, что граница всякого открытого множества, содержащего *x* и имеющего диаметр меньший, чем $ε_n$, состоит не менее чем из *n* точек [12].

Точка линии, индекс ветвления которой больше двух, называется точкой ветвления; индекс ветвления внутренних точек окружности равен двум; точка с индексом равным единице – концевая.

Следующие примеры значительно упрощают осмысление понятия. Пусть *O* – окружность единичного радиуса. Окружность *F* с центром в точке *G*, находящейся внутри *O*, пересекает ее в *R* = 2 точках, за исключением случаев, когда радиус *F* больше двух: при этом *R* = 0. Если *F* – граница некоторой окрестности точки *G*, не обязательно круглой, но «не слишком большой», то *R* равно, по меньшей мере, двум. Величина *R* = 2 называется степенью ветвления окружности, и остается неизменной для всех ее точек [4].

Линия, состоящая из отрезков, соединяющих некую точку со всеми точками канторова множества, лежащего в [0, 1], имеет во всех своих точках континуальный индекс ветвления [12].

Пусть множество *O* – треугольник Серпинского. В этом случае *R* не является одинаковой для всех точек *G*: во всех точках множества, за исключением вершин инициатора, значение *R* равно либо трем, либо четырем [4, 32].

Значение $R = 3$ характеризует любую точку множества $O$, являющуюся пределом бесконечной последовательности треугольников, каждый из которых содержится внутри треугольника предыдущего поколения и имеет вершины, отличные от вершин предшественника. Окружности, описанные вокруг этих треугольников, пересекают множество $O$ в трех точках, ограничивая при этом произвольно малые окрестности точки $G$. В этом случае, если $F$ ограничивает достаточно малую окрестность точки $G$, причем вершины инициатора должны лежать вне $F$, то, как показал В. Серпинский [32], $F$ пересекает $O$, как минимум в трех точках.

Значение $R = 4$ относится к вершинам любого конечного приближения к $O$ посредством треугольников. Вершина для аппроксимации порядка $h \geq k$ является общей вершиной $G$ для двух треугольников с длиной стороны $2^{-k}$. Окружности с центром в точке $G$ и радиусом $2^{-k}$ (при $h > k$) пересекают множество $O$ в 4 точках и ограничивают произвольно малые окрестности точки $G$. А если $F$ ограничивает «достаточно малую» окрестность точки $G$ (притом, что вершины инициатора также лежат вне $F$), то можно показать, что $F$ пересекает $O$, по меньшей мере, в 4 точках [4].

Когда множество $O$ – ковер Серпинского, результат оказывается радикально иным. Пересечение границы любой достаточно малой окрестности и $O$ представляет собой несчетное бесконечное множество точек, причем независимо от параметров $N$, $r$ и $D$. В этой дихотомии конечного/бесконечного треугольник Серпинского немногим отличается от стандартных кривых, а ковры Серпинского неотличимы от плоскости [4].

**Декартово произведение**

Прямое, или декартово, произведение – одна из основных общематематических конструкций, идея которой принадлежит

Декарту [45]. Декартовым произведением двух непустых множеств $X$ и $Y$ называется множество $X \times Y$, состоящее из всех упорядоченных пар вида $(x, y)$, где

$$x \in X, y \in Y : X \times Y = \{(x, y) \mid x \in X, y \in Y\}.$$

Если одно из множеств $X$ или $Y$ пусто, то произведение пусто. Множество $X \times Y$ можно отождествить с множеством функций, определенных на двухэлементном множестве $\{1, 2\}$ и принимающих значения в множестве $X$ при значении аргумента, равном 1, и в множестве $Y$ при значении аргумента, равном 2. Это отождествление позволяет распространить определение декартова произведения на случай любого количества множителей.

Значение конструкции прямого произведения определяется, прежде всего, тем, что в нем естественно вводится дополнительная структура, если все множители являются однотипными математическими структурами. Например, пусть $X_i$ – однотипные алгебраические системы, а именно, множества с общей сигнатурой конечноместных предикатов и операций, т.е. с общей совокупностью отношений и операций, действующих на основном множестве данной алгебраической системы, вместе с указанием их арностей, тогда произведение $X_i = \Pi X_i$ превращается в алгебраическую систему с той же сигнатурой. При этом выполнение во всех $X_i$ определенных тождеств влечет за собой их выполнение в произведении. Поэтому декартово произведение полугрупп, групп, колец, векторных пространств и т.п. снова являются полугруппами, группами, кольцами, векторными пространствами соответственно.

Многие задачи связаны с описанием математических объектов, неразложимых в декартово произведение, и с выяснением условий, при которых множители произведения определены однозначно с точностью до изоморфизма [45].

**Теория протекания**

Раздел теории вероятностей, имеющий собственные приложения в естественных и инженерных науках, – теория протекания, или перколяционная теория – изучает особенности возникновения и эволюции, а также свойства связных областей, т.н. бесконечных перколяционных кластеров [38, 46].

Перколяционная теория изучает континуальные и решётчатые задачи на плоскости и в объеме, а также $n$-мерные. Модели теории протекания просты и наглядны, некоторые задачи решаются аналитически, большинство – численно; часто используется метод Монте-Карло.

Простейшие решёточные задачи теории протекания формулируются следующим образом. Пусть есть решётка, рассматриваемая как совокупность узлов и связей. Каждый данный узел может быть помечен, например, чёрным цветом с вероятностью $x$. Совокупность связанных друг с другом чёрных узлов называется чёрным кластером. При $x = 0$ в системе нет чёрных кластеров, при $x \ll 1$ они представляют собой в основном совокупности малого количества узлов, а при $1 - x \ll 1$ в системе имеется чёрный бесконечный кластер (БК). Существует критическая концентрация $x_c$, при которой впервые возникает БК.

В компьютерной реализации перколяционных задач для определения места возникновения трансформированного участка – элемента будущего бесконечного кластера – используют генератор случайных чисел. Изменение концентрации таких участков в системе заметно модифицирует ее свойства. При достижении критической концентрации свойства системы изменяются скачкообразно: в матрице возникает перколяционный кластер, характерные размеры которого сравнимы с ее размерами, корреляционная длина расходится, меняется симметрия объекта, т.е. происходит структурный фазовый переход, т.н. переход второго рода [4, 17, 38, 46].

Параметром порядка перехода является мощность БК – веро-ятность того, что узел принадлежит бесконечному кластеру. Критическое поведение этой величины при $x \to x_c$ ($x > x_c$) определяется соотношением $P \propto (x - x_c)^\beta$, где $x_c$ – критическая концентрация конечных кластеров, $\beta$ – индекс параметра порядка. Характерный пространственный масштаб системы задаёт длина корреляции $\xi$. В первом приближении $\xi$ – это характерный размер конечных кластеров при $x < x_c$, и характерный размер пустот в БК при $x > x_c$. Критическое поведение этой величины определяется соотношением $\xi \propto |\tau|^{-\nu}$, где $\tau = (x - x_c)/x_c$, $\nu$ – индекс корреляционной длины [38].

Результаты фазового перехода могут разниться в зависимости от структуры и свойств бесконечного кластера [38].

## ЛИТЕРАТУРА

1. Mandelbrot B.B. Les objets fractals: Forme, hasard et dimension. – Paris: Flammarion, 1975. – 208 p.
2. Mandelbrot B.B. Fractals: Form, Chance, and Dimension. – San Francisco: W.H. Freeman and Company, 1977. – 352 p.
3. Mandelbrot B.B. The Fractal Geometry of Nature. – San Francisco: W.H. Freeman and Company, 1982. – 460 p.
4. Мандельброт Б. Фрактальная геометрия природы. – М.: ИКИ, 2002. – 656 с.
5. Александров П.С. Введение в общую теорию множеств и функций. – М.-Л.: ОГИЗ, 1948. – 413 с.
6. Александров П.С. Введение в теорию множеств и общую топологию. – М.: Наука, 1977. – 368 с.
7. Колмогоров А.Н., Фомин С.В. Элементы теории функций и функционального анализа. – М.: Наука, 1976. – 543 с.
8. Кроновер Р. Фракталы и хаос в динамических системах. Основы теории. – М.: Постмаркет, 2000. – 352 с.
9. Кузнецов С.П. Динамический хаос. – М.: Физматлит, 2006. – 356 с.
10. Лоскутов А.Ю. Динамический хаос. Системы классической механики. // УФН. – 2007. – Т. 177, вып. 9. – С. 989-1015.

11. Математическая энциклопедия в 5-ти томах. / Под ред. Ю.В. Прохорова. – Т. 2. – М.: СЭ, 1979. – 552 с.
12. Математическая энциклопедия в 5-ти томах. / Под ред. Ю.В. Прохорова. – Т. 3. – М.: СЭ, 1982. – 592 с.
13. Титчмарш Э. Теория функций. – М.: Наука, 1980. – 463 с.
14. Фихтенгольц Г.М. Курс дифференциального и интегрального исчисления, т. II. – М.: Наука, 1970. – 800 с.
15. Шредер М. Фракталы, хаос, степенные ряды. – Ижевск: НИЦ РХД, 2001. – 528 с.
16. Морозов А.Д. Введение в теорию фракталов. – Москва-Ижевск: ИКИ, 2002. – 160 с.
17. Федер Е. Фракталы. – М.: Мир, 1991. – 254 с.
18. Пуанкаре А. Наука и гипотеза. // Пуанкаре А. О науке – М.: Наука, 1990. – С. 38-78.
19. Герега А.Н. Об одном критерии относительной степени упорядоченности изображений. // Журнал технической физики. – 2010. – Т. 80, вып. 5. – С. 149-150.
20. Kullback S., Leibler R.A. On information and sufficiency. // Annals of Mathematical Statistics. – 1951. – V. 22, No. 1. – P. 79-86.
21. Кульбак С. Теория информации и статистика. – М.: Наука, 1967. – 480 с.
22. Адамар Ж. Элементарная геометрия. Планиметрия. – М.: Учпедгиз, 1948. – 609 с.
23. Корн Г., Корн Т. Справочник по математике для научных работников и инженеров. – М.: Наука, 1974. – 832 с.
24. Дирак П. Принципы квантовой механики. – М.: Наука, 1979. – 480 с.
25. Halsey T.C., Jensen M.N., Kadanoff L.P., Procaccia I., Shraiman B.I. Fractal measures and their singularities: The characterization of strange sets. // Phys.Rev. – 1986. – V. A33. – P. 1141-1151.
26. Hentschel H.G.E., Procaccia I. The infinite number of generalized dimensions of fractals and strange attractors. // Physica. – 1983. – V. 8. – 435-444.
27. Grassberger P. Generalized dimensions of strange attractors. // Phys. Lett. A. – 1983. – V. 97. – P. 227-230.
28. Гринченко В.Т., Мацыпура В.Т., Снарский А.А. Введение в нелинейную динамику. Хаос и фракталы. – М.: ЛКИ, 2010. –

280 с.
29. Божокин С.В., Паршин Д.А. Фракталы и мультифракталы. – Ижевск: НИЦ РХД, 2001. – 128 с.
30. Rényi A. Probability theory. – Amsterdam: North-Holland *Publ*. Co., 1970. – 666 p.
31. Юргенс Х., Пайтген Х.-О., Заупе Д. Язык фракталов. // В мире науки. – 1990. – № 10. – С. 36-44.
32. Sierpiński W. O krzywej, której każdy punkt jest punktem rozgałęzienia. // Prace matematyczno-fizyczne. – 1916. – Vol. 27. – P. 77-86.
33. Математическая энциклопедия в 5-ти томах. / Под ред. Ю.В. Прохорова. – Т. 5. – М.: СЭ, 1985. – 623 с.
34. Barnsley M.F., Rising H. Fractals Everywhere. – Boston: Academic Press, 1993. – 531 p.
35. Газале М. Гномон. М.: ИКИ, 2002. – 271 с.
36. Sierpinski W. Sur une courbe cantorienne qui contient une image biunivoque et continue de toute courbe donnee. // Comptes Rendus. – Series I. Mathematics. – 1916. – Vol. 162. – P. 629-632.
37. Урысон П.С. Труды по топологии и другим областям математики. Т. 2 / Ред., примеч. и вступ. ст. П.С. Александрова. – М.-Л.: ГИТТЛ, 1951. – С. 513-992.
38. Соколов И.М. Размерности и другие критические показатели в теории протекания. // УФН. – 1986. – Т. 150, вып. 2. – С. 221-255.
39. Ben-Avraham D., Havlin S., Movshovitz D. Infinitely ramified fractal lattices and percolation. // Philosophical Magazine B. – 1984. – V. 50, No. 2. – P. 297-306.
40. Герега А.Н., Дрик Н.Г., Угольников А.П. Ковер Серпинского с гибридной разветвленностью: перколяционный переход, критические показатели, силовое поле. // УФН. – 2012. – Т. 182, вып. 5. – С. 555-557.
41. Reynolds P.J, Stanley H.E., Klein W. Large-cell Monte Carlo renormalization group for percolation. // Physical Review B. – 1980. – V. 21, No. 3 – P. 1223-1245.
42. Menger K. Kurventheorie. – Leipzig: B.G. Teubner, 1932. – 375 s.
43. Halmos P.R. Naive Set Theory. – N.Y.: Springer, 1974. – 104 p.

44. Урысон П.С. Труды по топологии и другим областям математики. Т.1. / Ред., примеч. и вступ. ст. П.С. Александрова. – М.-Л.: ГИТТЛ, 1951. – 512 с.
45. Математическая энциклопедия в 5-ти томах. / Под ред. Ю.В. Прохорова. – Т. 4. – М.: СЭ, 1984. – 608 с.
46. Эфрос А.Л. Физика и геометрия беспорядка – М.: Мир, 1982 – 176 с.

## ПРИЛОЖЕНИЕ

## РАЗМЕРНОСТИ: ГЕНЕЗИС ПРЕДСТАВЛЕНИЙ И ФИЗИЧЕСКИЕ ПРИЛОЖЕНИЯ

### Введение

Теория размерности, основы которой, как известно, были заложены в работах А. Пуанкаре, А.Л. Лебега и Л. Брауэра [1-3], опубликованных в 1911-13 годах, открыла «доступ к изучению ряда интересных свойств точечных множеств, к построению обширной теории…» [4, 5].

Понятие размерности оказалось полезным и удобным концептом и в теоретической физике: во второй половине двадцатого века сформировались две основные группы задач, в которых размерность стала реальным инструментом получения и анализа решений. Первая – это задачи стохастической динамики, в частности, проблемы турбулентности. В этих задачах размерности являются доступными для измерения и структурно устойчивыми характеристиками системы, связанными, в частности, с показателями Ляпунова; позволяют провести классификацию странных аттракторов и связанного с ними хаотического поведения [6, 7].

Вторая группа – задачи перколяционной теории: раздела статистической физики, который на протяжении полувека изучает критические явления [6, 7]. Перколяционная теория адекватно описывает особенности возникновения и эволюции, а также свойства связных областей в системах, в которых имеет место геометрический фазовый переход. Она нашла применение в широком круге научно-технических задач: исследовании белковых структур, пористых тел, создании фильтров, изучении легированных полупроводников, при борьбе с эпидемиями, в исследованиях процессов полимеризации, при создании композиционных материалов, изучении мировоззренческих вопросов и мно-гих других. При этом для критических показателей физических величин, описывающих процессы и явления, как правило, можно указать множество, с размерностью которого этот показатель связан. В свою очередь, исследование структуры этих множеств много даёт для понимания критического поведения системы и соотношений между показателями, позволяет проследить связь между поведением системы в промежуточной асимптотике и её геометрией [6, 7].

Перколяционная теория известна также как раздел теории вероятностей, имеющий собственные приложения в естественных и инженерных науках [8-11].

В математических и физических исследованиях используется большое количество размерностей, и в конкретном исследовании всегда встаёт вопрос выбора, который обусловлен содержанием задачи. Размерности, представленные в обзоре, соотносятся с совершенно разными понятиями, в частности, имеют отношение и к описанным группам задач (в первую очередь, к перколяционной теории), и к топологическим исследованиям.

Далее в тексте, как и во всех физических приложениях, размерности определяются как показатель степени в выражениях типа $a \sim b^c$.

**Геометрическая и физическая размерности**

Размерность физической величины определяет её связь с величинами, положенными в основу системы единиц измерения, т. е. устанавливает соотношение масштабов данной и основных единиц измерения [12]. Из такого определения следует, что физические величины, в отличие от геометрических, изменяются не только при преобразованиях координат, но и при модификации системы физических величин. Это значит, что при определённом выборе основных единиц каждому физическому объекту соответствует один и только один геометрический объект – образ физического объекта, который приобретает множитель $a^{\alpha} \beta^{b}$ $c^\gamma \ldots l^\lambda\, m^\mu\, n^\nu$, если вводится другое множество допустимых основных единиц; этот множитель называется размерностью геометрического образа, или абсолютной размерностью физического объекта. Абсолютная размерность определена относительно аффинной группы, если мы работаем с прямолинейными координатами $n$-мерного евклидового пространства. Эта размерность не является той размерностью, с которой работают физики [13].

Если компоненты физического объекта относительно локальной декартовой системы координат, основанной на единице длины, приобретают при изменении основных единиц множитель $a^{\alpha 1}\, \beta^{b1}\, c^{\gamma 1} \ldots l^{\lambda 1}\, m^{\mu 1}\, n^{\nu 1}$, то его называют относительной размерностью, или просто размерностью, физического объекта. Это та размерность, которая используется в физике.

Таким образом, различие между абсолютной (геометрической) и относительной (физической) размерностью физического объекта обусловлено тем, что общая система координат не связана с единицей длины, а локальная декартова изменяется при введении другой единицы. Следовательно, соотношение между ними зависит только от закона трансформации объекта при преобразовании координат [13].

**Размерность пространства. Топологическая размерность**

Если рассматривать геометрические объекты как множества точек евклидова пространства $R^E$, то понятие топологической размерности можно ввести по рекуррентной схеме, предложенной А. Пуанкаре [1].

Положим размерность любого конечного или счётного множества точек равным нулю; размерность связного множества будем считать равной $d + 1$, если оно может быть разрезано на две несвязные части исключением из него как минимум $d$-мерного множества точек, т. е. проведением $d$-мерного разреза. При таком определении, если положить размерность точки рав-ной нулю, то топологическая размерность линии будет равна единице, плоскости и сферы – двум, шара – трём и т. д.

Из определения видно, что топологическая размерность может быть только целым числом, и совпадает с интуитивным представлением о минимальном количестве переменных, которые нужно задать для определения положения точки на объекте [14]. В 1902 году в книге «Наука и гипотеза» [15] А. Пуанкаре писал, что «размерность пространства – это минимальное число параметров, которые необходимы, чтобы отличать точки пространства друг от друга»; полное зрительное пространство «имеет как раз три измерения; т. е. элементы наших зрительных ощущений … будут вполне определены, когда известны три из них». И если в пространстве это число равно трём, на плоскости достаточно двух координат, на линии – одной; в этом смысле пространство – трёхмерно, плоскость – двумерна, линия – одномерна.

В статье «Почему пространство имеет три измерения» [1] Анри Пуанкаре определил размерность, которую можно было бы назвать антропной. Он пишет о существовании «экспериментальных фактов, которые заставляют нас приписывать пространству три измерения. Именно ввиду этих данных нам было удобнее приписать ему три измерения, а не четыре или два. Но слово «удобный», пожалуй, в данном случае недостаточно сильно: существо, которое приписывало бы пространству два или четыре измерения, оказалось бы менее приспособленным к борьбе за существование в мире, подобном нашему». В случае двух измерений оно предполагало бы существование таких соотношений, которые мы, люди, не допускаем; а в случае четырёх – отбрасывало бы такие, которыми мы пользуемся [1].

Физическое обоснование трёхмерности пространства дано в работе П. Эренфеста [16], в которой исследуется, в частности, аналог гравитационного закона Ньютона для пространств с различным числом измерений. В этом случае зависимость гравита-ционных сил от расстояния определяется выражением $F = G_i M m / R^{i-1}$, где $M, m, R$ – массы и расстояние в «классическом» понимании, $i$ – число пространственных координат, $G_i$ – коэффициенты, в частности, $G_3$ – гравитационная постоянная Ньютона. Некоторым основанием для такого предположения служат результаты анализа законов движения, полученные Эренфестом для пространств с числом измерений, отличным от трёх [16, 17].

В предложенном Ньютоном гравитационном законе сила пропорциональна $R^{-2}$, причём эта зависимость неоднократно проверялась экспериментально, и значение показателя степени установлено с точностью до $2 \pm 3 \cdot 10^{-11}$, следовательно, с этой же точностью размерность нашего пространства $i = 3$. Эти данные получены из прецизионных измерений орбиты Луны, движущейся вокруг Земли, и хотя средний радиус лунной орбиты равен 384 тысячам километров, модельные данные отличаются от измеренных на 4 мм [18].

**Адекватность меры. Размерность Хаусдорфа-Безиковича**

Для классических геометрических объектов понятие топологической размерности решает вопрос об адекватной мере.

Известный способ измерить длину кривых, площадь поверхностей и объёмы тел состоит в разделении пространства на малые кубы или сферы с характерным линейным размером $\delta$. Подсчитывая количество отрезков, квадратов или кубов, необходимых для покрытия рассматриваемого множества точек, можно получить меру этого множества [19].

Для обычной кривой длина $L$ может быть определена предельным переходом

$$L = N(\delta)\, \delta \underset{\delta \to 0}{\longrightarrow} L_0\, \delta^0,$$

где $\delta$ – длина прямолинейных отрезков, $N(\delta)$ – их количество. Как видно, в пределе при $\delta \to 0$ мера $L$ становится асимптотиче-ски равной длине кривой и не зависит от $\delta$.

Длина является адекватной мерой обычной (нефрактальной) кривой. Если поставить в соответствие линии не длину, а площадь или объём, то та же процедура покажет, что такие меры обращаются в нуль. Действительно, пусть $N(\delta)$ – количество квадратов, необходимых для покрытия кривой, $\delta^2$ – площадь одного квадрата, тогда

$$S = N(\delta)\, \delta^2 \xrightarrow[\delta \to 0]{} L_0\, \delta^1 = 0;$$

аналогично,

$$V = N(\delta)\, \delta^3 \xrightarrow[\delta \to 0]{} L_0\, \delta^2 = 0.$$

Центральное место в определении размерности Хаусдорфа-Безиковича занимает понятие адекватности меры.

Рассмотрим множество точек, образующих поверхность в трёхмерном пространстве. Адекватная мера такого множества – площадь. Действительно,

$$S = N(\delta)\, \delta^2 \xrightarrow[\delta \to 0]{} S_0\, \delta^0,$$

где $S_0$ – площадь поверхности. Таким образом, количество квадратов, необходимых для покрытия поверхности, определяется в пределе при $\delta \to 0$ как $N(\delta) = S_0/\delta^2$.

Если поставить в соответствие поверхности длину, то

$$L = N(\delta)\, \delta \xrightarrow[\delta \to 0]{} S_0\, \delta^{-1} = \infty,$$

что говорит о невозможности покрыть поверхность конечным количеством отрезков прямой. Если сделать попытку установить соответствие между поверхностью и объёмом, то он обратится в нуль

$$V = N(\delta)\, \delta^3 \xrightarrow[\delta \to 0]{} S_0\, \delta^1 = 0.$$

Определим меру фрактального множества, используя проб-ную функцию $M_d$, которая в зависимости от выбора её размерности $d$, обращается в нуль или бесконечность при $\delta \to 0$. Введём размерность Хаусдорфа-Безиковича $D_H$, при которой мера $M_d$ изменяет значение с нуля на бесконечность

$$M_d = \sum \gamma(d) \delta^d = \gamma(d) N(\delta) \delta^d \underset{\delta \to 0}{\to} \begin{cases} 0, & \text{при } d > D, \\ \infty, & \text{при } d < D, \end{cases}$$

где $\gamma(d)$ – геометрический коэффициент, зависящий от формы элементов, покрывающих множество; для квадратов и кубов он равен единице, для кругов – $\pi/4$, для сфер – $\pi/6$ [15]. Для физиков такое поведение меры означает, что $D$ представляет собой критическую размерность [20].

Существенно, что при определении размерности Хаусдорфа-Безиковича необходимо покрывать множество элементами всевозможных размеров, не превышающих некоторое малое значение, и определить infimum выражения $\gamma(d) \sum \delta^d$. Очевидно, что процесс минимизации этой суммы по всем возможным разбиениям чрезвычайно трудоёмок, и обычно производят оценку размерности Хаусдорфа-Безиковича величиной ёмкости множества $D_c$. Это типичная ситуация в прикладных задачах теории размерностей: среди однотипных иногда можно найти размерности, пригодные для оценки значений других, расчёт которых трудоёмок или нереализуем.

Для определения $D_c$ рассмотрим случай покрытия множества точек в $d$-мерном евклидовом пространстве минимальным количеством $d$-мерных кубиков (или сфер) одинакового размера. (Покрытие сферами используется для того, чтобы не говорить об ориентации). То есть, если $N(\delta) \sim \delta^{-D_c}$, то $D_c$ – колмогоровская ёмкость множества [6, 21]; индекс $c$ – сокращение от англ. capacity.

Пусть $A$ есть величина, характеризующая покрытие, и $N(\delta) \approx A \cdot \delta^{-d}$ – минимальное количество $d$-мерных кубиков. Логарифмируя это выражение, получим

$$ln\ N(\delta) \approx ln\ A - d \cdot ln\ \delta,$$

откуда, приблизительно,

$$d = -\frac{ln\ N(\delta)}{ln\ \delta} + \frac{ln\ A}{ln\ \delta}\ .$$

Так как $ln\ \delta \to -\infty$ при $\delta \to +0$, то ёмкость множества есть предел

$$D_c = -\lim_{\delta \to 0}\frac{ln\ N(\delta)}{ln\ \delta},$$

который, обычно, существует.

Поскольку при определении хаусдорфовой размерности должны использоваться всевозможные покрытия множества, а при расчёте ёмкости – элементы одного размера, то $D_H < D_c$.

Хаусдорфова размерность и ёмкость по Колмогорову могут различаться даже для очень простых множеств [6]. Например, для множества точек прямой с координатами $x_n = 1/n$ первая равна 0, вторая – 1/2. Ёмкости, в отличие от размерностей, в частности, хаусдорфовой, не остаются инвариантными при кусочно-гладком, возможно, имеющем особенности, преобразовании координат, а для величины, определяемой как размерность, такая инвариантность необходима.

Как отмечалось, во всех физических приложениях размерность определяется как показатель $M \sim l^D$, где $M$ – некое свойство, $l$ – характерный размер, а определить является ли она хаусдорфовой размерностью или ёмкостью не представляется возможным. Это связано с тем, что размерность описывает свойства промежуточной асимптотики, и переход к пределу, требуемый формальными определениями, невозможен. Кроме того, на малых масштабах система не является фрактальной, и её поведение описывается некоторым минимальным масштабом.

Так определённые размерность Хаусдорфа-Безиковича и ёмкость есть локальные свойства в том смысле, что характеризуют множество точек при исчезающе малом характерном размере пробной функции ($\delta \to 0$), используемой для его покрытия.

Важно, что для простых геометрических объектов хаусдорфова размерность совпадает с топологической. Действительно, пусть есть квадрат со стороной $a$; покроем его малыми квадратами площадью $\delta^2$. Тогда для меры $A_d(\delta)$ получим

$$A_d(\delta) = \sum \delta_i^d = N(\delta) \cdot \delta^d \approx a^2 \delta^{d-2}.$$

При $d < 2$ мера $A_d(\delta)$ неограниченно возрастает при уменьшении $\delta$, в случае $d > 2$ – стремится к нулю. Следовательно, по определению размерности Хаусдорфа-Безиковича $D_H = 2$.

**Размерность Хаусдорфа-Безиковича как фрактальная**

В евклидовом пространстве $R^E$ величина топологической размерности $D_T$ и размерности Хаусдорфа-Безиковича $D_H$ заключены в промежутке между нулем и $E$. При этом топологическая всегда является целым числом, а для размерности Хаусдорфа-Безиковича это не обязательно. Для евклидовых множеств $D_H = D_T$, в общем же случае эти две размерности удовлетворяют неравенству Шпилрайна (Edward Szpilrajn) $D_H \leq D_T$ [5, 20].

Однако существуют множества, для которых $D_H > D_T$. В [20] Б. Мандельброт пишет: «Такие множества необходимо было как-то называть, поэтому я придумал термин «фрактал», определив его следующим образом: фракталом называется множество, размерность Хаусдорфа-Безиковича для которого строго больше его топологической размерности».

Любое множество с нецелым значением $D_H$ является фракталом; фрактал может иметь и целочисленную размерность. Если понимать термин «дробь» как синоним выражения «нецелое вещественное число», то часто значения размерности $D_H$ являются дробными. Учитывая, что $D_H$ может принимать и целочисленные значения $D_T < D_H \leq E$, Бенуа Мандельброт предпочёл назвать её фрактальной размерностью [20], и обозначить через *D*.

**Многообразие покрытий. Размерности Минковского-Булигана и Понтрягина-Шнирельмана**

По Мандельброту фрактальная размерность и все её возможные варианты – не топологические, но метрические понятия: каждая включает в себя метрическое пространство, в котором определены расстояния между любыми двумя точками [20]. При этом сами размерности определяются алгоритмами покрытия множества *d*-мерными шарами и, в сущности, есть функции способа покрытия.

В способе покрытия ограниченного множества, предложенном Г. Кантором, каждая его точка рассматривается как центр шара [20]. Такой подход связан с очевидными неудобствами. Во-первых, для множеств, содержащих бесконечное количество точек, такой алгоритм неоперабелен. Ситуация, однако, разрешается тем, что, оказывается, достаточно построить конечное число шаров $N(\rho) \sim 1/\rho$.

Во-вторых, сумма перекрывающихся объёмов шаров при $\rho \to 0$ не должна непременно сходиться к протяжённости (в смысле Минковского) множества, т. е. *d*-мерному объёму: длине, площади, объёму и т. д. В примере Х. А. Шварца [22] показано, что по мере увеличения точности триангуляции боковой поверхности прямого кругового цилиндра, сумма площадей треугольников не обязательно сходится к её площади, а может быть равной сколь угодно большой конечной или бесконечной величине.

При рассмотрении этого парадокса Г. Минковский показал, что если определить протяжённость как

$$V\{d\text{-мерный шар радиуса } \rho\} = \gamma(d)\,\rho^D,$$

где *d* – стандартная топологическая размерность рассматриваемого множества, множитель $\gamma(d) = [\Gamma(½)]^d / \Gamma(1 + d/2)$, $\rho$ – ра-диус покрывающих шаров, $D$ – размерность фрактального множества (если множество нефрактально, $D = d$), то при $\rho \to 0$ она может не иметь предела [20]. В этом случае выражение $\lim_{\rho \to 0} V$ заменяется на $\lim_{\rho \to 0} \sup V$ и $\lim_{\rho \to 0} \inf V$ – верхнюю и нижнюю протяжённости множества. При этом любому вещественному числу из интервала $]\liminf, \limsup[$ соответствует, по меньшей мере, одна последовательность значений $\rho_m \to 0$, таких, что при $m \to \infty$ сумма площадей треугольников в примере Х. А. Шварца сходится к площади поверхности [20]. Г. Минковский показал также, что в случае стандартных евклидовых структур существует величина $D_M$ – размерность Минковского – такая, что при $d > D_M$ верхняя протяжённость множества обращается в нуль, а при $d < D_M$ нижняя – бесконечна [20].

В 1928 году Ж. Булиган обобщил размерность Минковского на случай дробных $d$, и показал, что она определяется выражением $\lim_{\rho \to 0} \inf V$. В некоторых случаях величины размерностей Минковского-Булигана $D_{MB}$ и Хаусдорфа-Безиковича $D_H$ совпадают, например, для гладких кривых и поверхностей, и с учётом того, что $D_{MB}$ легче поддаётся оценке, и по аналогии с ситуацией с размерностью Хаусдорфа-Безиковича и ёмкостью, она может использоваться для определения величины $D_H$ [20].

Несложно привести и обратный пример, когда такая оценка невозможна: для компактного множества $\{0, 1, 1/2, 1/3, 1/4, \ldots\}$ как для всякого счётного $D_H = 0$, а $D_{MB} = 1/2$. В общем случае, как показано в [23], $D_{MB} \geq D_H$.

Среди всевозможных наборов покрывающих шаров наиболее экономичным является комплект, содержащий минимум шаров $N(\rho)$, который используют для определения размерности Понтрягина-Шнирельмана [20, 24]

$$D_{PSch} = \lim_{\rho \to 0} \inf \ln N(\rho) / \ln(1/\rho).$$

**Размерность самоподобия и клеточная размерность**

Как известно, объекты инвариантные относительно изменения масштаба и параллельного переноса называются самоподобными. Если при соответствующем изменении масштаба в $n$ раз ($n < 1$) можно однократно покрыть исходный объект уменьшенными копиями, то он самоподобен с коэффициентом подобия $n$ $(N) = 1/N^{1/d}$, где $N$ – количество одинаковых частей, имеющих в $n$ раз меньший линейный размер, $d$ – размерность подобия (самоподобия), равная топологической размерностью объекта. В случае геометрически самоподобных (регулярных) фракталов

$$n\,(N) = 1/N^{1/D_S},$$

где $D_S$ совпадает с размерностью Хаусдорфа-Безиковича [19] и определяется формулой

$$D_S = \ln N / \ln n.$$

В качестве иллюстрации рассмотрим квадрат, разделённый на $N = 4$ равных квадрата со сторонами в $n = 2$ раза меньшими, чем у исходного. Тогда размерность самоподобия, равная топологической, для квадрата имеет значение 2; для куба, разделённого на $N = 8$ равных частей, $D_S = \ln 8 / \ln 2 = 3$. Для кривой Коха [19] – регулярного фрактала, при построении которого на каждой итерации масштаб покрывающих отрезков уменьшается в $n = 3$ раза, а их количество становится равным $N = 4$ размерность $D_S = \ln 4 / \ln 3 = 1{,}2618\ldots = D_H$.

Для определения размерности нерегулярных объектов фрактального типа, например, изображений государственной границы или береговой линии, описанный алгоритм, естественно, не подходит, и применяют другое определение размерности, связанное с иным алгоритмом.

Пространство, в котором расположен интересующий нас объект, разбивают на клетки размером $\delta^2$, например, наносят с помощью палетки на изображение объекта квадратную сетку со стороной $\delta$, и подсчитывают число клеток, содержащих точки объекта. Разбиение многократно повторяют, используя всё меньший масштаб. Зависимость количества клеток, в которые попали точки объекта, от $\delta$ описывается выражением $N(\delta) = A\delta^{-D_S}$, где $D_S$ – искомая фрактальная размерность самоподобия. Для расчёта её значения строят график зависимости $N(\delta)$ в двойном логарифмическом масштабе, при этом угловой коэффициент графика определяет значение $D_S$.

За размерностью, определённой по такому алгоритму, можно сохранить название размерности самоподобия: действительно, если применить описанный алгоритм к регулярным объектам, то значения размерностей совпадают с вычисленными по формуле для $D_S$. Однако как размерность, определяемую посредством подсчёта количества клеток, её называют клеточной [6, 19].

**Размерность энтропии меры**

Описанный алгоритм расчёта размерности нерегулярного фрактального множества имеет естественное ограничение: если, например, береговая линия сильно изрезана и неоднократно пересекает некую клетку, то в количество клеток, покрывающих множество точек, она всё равно даёт единичный вклад, что «не вполне честно» [19].

Рассмотрим распределение точек множества по клеткам, отражающее распределение меры. Пусть множество, состоящее из $N$ точек, имеет в $i$-й клетке $N_i$ точек. И пусть $\mu = N_i / N$ – вероятность заполнения клетки. Можно построить меру

$$M_d(q, \delta) = \sum_{i=1}^{N} \mu_i^q \delta^d = N(q, \delta)\delta^d \xrightarrow[\delta \to 0]{} \begin{cases} 0, & \text{при } d > \tau(q) \\ \infty, & \text{при } d < \tau(q) \end{cases}$$

обладающую показателем $d = \tau(q)$, при котором она не обращается в нуль или бесконечность при $\delta \to 0$ [19]. Мера характеризуется всей последовательностью показателей $\tau(q)$, определяющих степенной закон, по которому изменяются вероятности $\{\mu\}$ в зависимости от $\delta$. При этом взвешенное число клеток равно

$$N(q,\delta) = \sum_i \mu_i^q \sim \delta^{-\tau(q)},$$

где

$$\tau(q) = -\lim_{\delta \to 0} \frac{\ln N(q,\delta)}{\ln \delta}.$$

Из этих соотношений видно, что при $q = 0$ получаем $\mu_i^{q=0} = 1$. Тогда $N(q=0,\delta)$ – количество точек, покрывающих множество, и $\tau(0) = D$ – фрактальная размерность множества. Кроме того, с учётом $\sum_i \mu_i = 1$ получаем $\tau(1) = 0$.

Введём производную

$$\frac{d\tau(q)}{dq} = -\lim_{\delta \to 0} \frac{\sum_i \mu_i^q \ln \mu_i}{(\sum_i \mu_i^q) \ln \delta},$$

и рассмотрим

$$\left.\frac{d\tau(q)}{dq}\right|_{q=1} = -\lim_{\delta \to 0} \frac{\sum_i \mu_i \ln \mu_i}{\ln \delta} = \lim_{\delta \to 0} \frac{S(\delta)}{\ln \delta},$$

где $S(\delta)$ – информационная энтропия разбиения меры $M$ по ячейкам размера $\delta$, которую можно записать в виде

$$S(\delta) = -\sum_i \mu_i \ln \mu_i \sim -\alpha_1 \ln \delta.$$

Показатель $\alpha_1 = -\left.\frac{d\tau(q)}{dq}\right|_{q=1}$ есть фрактальная размерность множества, на котором сосредоточена мера; он описывает скейлинговое поведение энтропии разбиения меры при изменении размера ячейки $\delta$ [19]; с точностью до множителя она равна информационной размерности – второй из спектра обобщённых размерностей Реньи.

**Обобщённые размерности Реньи**

Для адекватного описания неоднородных фрактальных объ-ектов (мультифракталов) можно использовать спектр обобщённых фрактальных размерностей А. Реньи [14, 25, 26].

Пусть есть фрактальный объект, ограниченный произвольной областью размера $L$ в евклидовом пространстве размерности $d$. И пусть он представляет собой множество $N \to \infty$ точек, произвольно распределённых в этой области. Разобьём всю область на прямоугольные ячейки со стороной $\delta \ll L$ и объёмом $\delta^d$. Ячейки, в которых содержится хотя бы одна точка, определим как занятые. Пусть $N(\delta)$ – суммарное количество заполненных ячеек, $n_i(\delta)$ – количество точек в $i$-й ячейке; определим вероятность того, что произвольная точка множества находится в ячейке $i$ как

$$p_i(\delta) = \lim_{N \to \infty} (n_i(\delta)/N),$$

и введём обобщённую статистическую сумму

$$Z(q,\delta) = \sum_{i=1}^{N(\delta)} p_i^{\,q}(\delta).$$

Тогда по А. Реньи [25, 26] можно ввести спектр обобщённых размерностей, характеризующих распределение точек в произвольной области, и показывающих насколько оно неоднородно. Существенно, что эти размерности связаны с показателями $\tau(q)$, описанными в предыдущем разделе, как

$$D_q = \frac{\tau(q)}{q-1},$$

где $\tau(q) = \lim_{\delta \to 0}(\ln Z(q,\delta)/\ln \delta)$, а множитель $1/(q-1)$ выбран, чтобы для множеств постоянной плотности в $E$-мерном пространстве выполнялось равенство $D_q = E$ [19].

Действительно, если для равномерно распределённой меры в $E$-мерном пространстве с постоянной плотностью точек разделим пространство на $N = \delta^{-E}$ ячеек объёмом $\delta^E$, тогда $\mu_i = \delta^E$ и

$$\sum_{i=1}^{N} \mu_i^q = \sum_{i=1}^{N} \delta^{qE} = \delta^{(q-1)E},$$

и, следовательно,

$$D_q = \frac{1}{q-1} \lim_{\delta \to 0} \frac{\ln \delta^{(q-1)E}}{\ln \delta} = E.$$

Таким образом, спектр фрактальных размерностей для равномерно распределённой меры сводится к размерности пространства и не зависит от порядка момента $q$ [19].

**Информационная и корреляционная размерности. Свойства функции $D_q$**

Определим смысл обобщённых размерностей Реньи $D_q$ для $q = 1$ и $q = 2$ [14, 25, 26].

Обобщённая статистическая сумма в силу условия нормировки вероятности при $q = 1$ равна единице, что, очевидно, приводит к неопределённости 0 / 0 в выражении для $D_1$. Раскроем неопределённость с помощью выражения

$$Z(q,\delta) = \sum_{i=1}^{N(\delta)} p_i^q = \sum_{i=1}^{N(\delta)} p_i \, exp[(q-1)\ln p_i].$$

Устремляя $q \to \infty$ и раскладывая в ряд экспоненту, с учётом условия нормировки, получим

$$Z(q \to 1, \delta) \approx \sum_{i=1}^{N(\delta)} [p_i + (q-1)p_i \ln p_i] = 1 + (q-1)\sum_{i=1}^{N(\delta)} p_i \ln p_i.$$

Тогда

$$D_1 = \lim_{\delta \to 0} \sum_{i=1}^{N(\delta)} p_i \ln p_i / \ln \delta.$$

С точностью до знака числитель в этой формуле представляет собой энтропию $S(\delta)$ фрактального множества

$$S(\delta) = -\sum_{i=1}^{N(\delta)} p_i \ln p_i,$$

и, следовательно, формула для $D_1$ приобретает вид

$$D_1 = -\lim_{\delta \to 0} \frac{S(\delta)}{ln(\delta)}.$$

Таким образом, величина $D_1$ характеризует информацию, необходимую для определения положения точки в некоторой ячейке, и показывает, как возрастает количество информации при стремлении размера ячейки к нулю $S(\delta) \approx \delta^{-D_1}$. Её называют информационной размерностью [14, 25].

Определение $D_1$ через энтропию даёт возможность сравнить величины фрактальной и информационной размерностей: расчёт $D_1$ имеет смысл лишь в случае неоднородности фрактального множества точек, для которого энтропия меньше, чем для однородного, и, следовательно, всегда $D_1 < D_0$ (равенство, очевидно, возможно лишь для однородного случая). В теории размерности этот результат обобщён для произвольного $q$, и показано, что имеет место неравенство

$$D_0 \geq D_1 \geq D_2 \geq D_3 \ldots$$

Для определения физического смысла обобщённой размерности $D_2$, равной

$$D_2 = \lim_{\delta \to 0} \left( ln \sum_{i=1}^{N(\delta)} p_i^2 / ln(\delta) \right),$$

введём парный корреляционный интеграл

$$I(\delta) = \lim_{N \to \infty} \frac{1}{N^2} \sum_{n,m} \theta(\delta - |r_n - r_m|).$$

Здесь суммирование проводится по всем парам точек фрактального множества с радиус-векторами $r_n$ и $r_m$, а $\theta(x)$ – ступенчатая функция Хэвисайда. Эта сумма определяет количество пар точек $(n, m)$, расположенных на расстоянии меньшем, чем $\delta$, и, следовательно, будучи разделённой на $N^2$, равна вероятности того, что две произвольные точки обладают этим свойством.

Эту же вероятность можно найти иначе [14, 25, 26]. По определению величина $p_i^2$ – вероятность попадания двух точек в $i$-ю ячейку с размером $\delta$. Суммируя $p_i^2$ по всем занятым ячейкам, получим вероятность того, что две произвольно выбранные точки из множества лежат внутри одной ячейки с размером $\delta$. Следовательно, расстояние между ними будет порядка $\delta$ или меньше. Тогда, с точностью до численных коэффициентов

$$I(\delta) \approx \sum_{i=1}^{N(\delta)} p_i^2 \approx \delta^{D_2}.$$

Таким образом, обобщённая размерность $D_2$ определяет зависимость корреляционного интеграла $I(\delta)$ в пределе $\delta \to 0$, и называется корреляционной [25].

С помощью вероятностной интерпретации удобно выяснить смысл граничных размерностей спектра Реньи: минимальной $D_\infty = \lim_{q \to \infty} D_q$ и максимальной $D_{-\infty} = \lim_{q \to -\infty} D_q$.

Величина $D_{-\infty}$ определяет формальную верхнюю границу интервала изменений $D_q$ (максимальное значение имеет размерность $D_0$, т. к. размерности в спектре определены для неотрицательных $q$), размерность $D_\infty$ – минимальная из размерностей. При $q \to \infty$ основной вклад в обобщённую статистическую сумму, очевидно, вносят ячейки с наибольшей вероятностью заполнения, т.е. содержащие максимальное число частиц; при $q \to -\infty$, соответственно, – самые разрежённые. Тогда, с учётом описанного выше,

$$D_{-\infty} \geq D_0 \geq D_1 \geq D_2 \geq D_3 \geq \ldots \geq D_\infty.$$

**Внешняя и внутренняя размерности кривой**

Если применить идеи Хаусдорфа для определения размерности фрактальной кривой $\gamma$, и для этого строить вокруг её точек кружки радиуса $\delta \to 0$ и вычислять площадь их объединения $S(\delta)$ (учитывая площадь перекрытия нескольких кружков только один раз), то скорость убывания площади с уменьшением $\delta$ определит размерность кривой [28].

Действительно, для гладкой кривой $S(\delta) \sim \delta \cdot L$, где $L$ – её дли-на, для плоской области $S(\delta) \sim \delta^0$, а для фрактала оценка площади, которая определяется шириной окрестности кривой, зависит в этом случае не от радиуса кружка, а от размера изгибов кривой. Определим длину «волнового вектора» $k = 1/\delta$, тогда отдельные «периоды» изгибов кривой укладываются в кружке размера $\delta$, и если $a(k)$ при $\delta \to 0$ убывает медленнее, чем $\delta$, то ширина полосы кружков не порядка $\delta$, а порядка $a(1/\delta)$; при $a(1/\delta) < \delta$ полоса кружков успевает отслеживать все изгибы кривой, и она является не фрактальной, а гладкой [28].

Пусть кривая $\gamma$ – фрактальна. Тогда суммарная площадь объединений всех кружков $S(\delta) \sim a(1/\delta) \sim \delta^\alpha$, и если $0 < \alpha < 1$, то $S(\delta)$ убывает медленнее, чем для гладкой кривой, следовательно, кривая $\gamma$ занимает промежуточное положение между линией и площадью. Ф. Хаусдорф предложил определение, согласно которому размерность такого образования равна

$$dim_{ext} \gamma = 2 - \alpha.$$

Индекс *ext* (сокращение от *external* – англ. «внешний; наружный») указывает, что при построении этой величины нам пришлось выйти за пределы самой кривой [28].

По аналогии с двумерным случаем внешняя размерность фрактальной кривой в пространстве равна

$$dim_{ext} \gamma = 3 - 2\alpha,$$

в $n$-мерном –

$$dim_{ext} \gamma = n - (n-1)\alpha.$$

Определим внутреннюю размерность кривой $\gamma$. Для этого разделим её на участки длиной $\delta$ и введём параметр $t$. Вычислим длину кривой, учитывая лишь те её изгибы, на которых $t$ изменяется не менее чем на $\delta$. Получим, что сумма длин этих отрезков порядка $a(1/\delta)(1/\delta) \sim \delta^{\alpha-1}$ и стремится к бесконечности с уменьшением $\delta$ [28].

С чем связана эта расходимость. Предположим, вслед за авторами [28], что мы ошиблись в определении размерности на-шего объекта и исследуем не кривую, а пытаемся определить одним параметром плоскость. Такая параметризация, конечно, плоха – линия всё более плотно и с самопересечениями заполняет плоскость, образуя на ней подобие решётки. Расстояние между её полосами ~ $\delta$, а число квадратов ~$1/\delta^2$. Размеры звеньев ломаной, плотно устилающей плоскость, малы, но число их очень велико: сумма длин отрезков ломаной ~ $\delta \cdot 1/\delta^2 \to \infty$. И это естественно: так как область двумерна, нужно подсчитывать не её длину, а площадь, т.е. суммировать не длины сторон, а квадраты длин, сумма которых конечна [28].

В случае дробной размерности нужно суммировать некоторые $\mu$-е степени длин. Для конечности получающихся сумм нужно положить $\mu = 1/\alpha$. Размерность этой суммы равна $см^{1/\alpha}$, а само число $1/\alpha$ является размерностью. Тогда в качестве внутренней размерности фрактала естественно принять число

$$dim_{int} \gamma = 1/\alpha.$$

Эта формула сохраняется и для кривых в пространстве любого числа измерений.

Внешняя размерность кривой фрактального типа на плоскости изменяется от 1 до 2 (размерность пространства), а внутренняя – от 1 до бесконечности, и совпадают они только для тривиального случая гладкой кривой. В общем случае внешняя размерность фрактальной кривой изменяется от размерности гладкого объекта до размерности пространства, а внутренняя – от размерности гладкого объекта до бесконечности [28].

«Очевидно, что в разных физических задачах нужно пользоваться разными определениями фрактальной размерности. Например, если мы интересуемся задачей адсорбции на тонкую нитку, то для нас важно знать, сколько атомов сможет поместиться вблизи нитки, т.е. внешнюю размерность. Если же мы хотим оценить вес нитки, то важна размерность внутренняя» [28].

### Массовые размерности

В структуре вещества всегда можно выделить масштаб, равный корреляционной длине $\xi$, т. е. расстоянию, вне которого частицы вещества ведут себя статистически независимо, и которое определяет верхнюю границу промежуточной асимптотики, а значит, и границу между интенсивным и экстенсивным поведением плотности вещества тела. Действительно, при $l < \xi$ масса тела определяется соотношением $M \sim l^D$, где $D$ – фрактальная размерность, и выражением $M \sim l^d$ при $l > \xi$ (здесь $d$ – размерность пространства). Тогда плотность можно определить как

$$\rho = \begin{cases} M/l^d \sim l^{D-d}, & \text{при } l \leq \xi \\ const, & \text{при } l > \xi. \end{cases}$$

Воспользуемся этим соотношением, и рассчитаем, для примера, массу кубика сахара рафинада с ребром $L$. Существует альтернатива: определить массу через среднюю плотность $\rho$ как

$$m = \rho L^3,$$

или учесть, что песчинки сахара образуют статистически самоподобную структуру, и тогда масса равна

$$m = \rho_c L^{D_m},$$

где $\rho_c$ – плотность сахара, $D_m$ – массовая размерность.

В отличие от «сплошного» тела во фрактальном объекте средняя плотность зависит от объёма, т.е. является экстенсивной физической величиной, убывающей при его возрастании. Последнее обстоятельство легко объясняется редко используемым, «ненаучным» определением фрактала, предложенным Б. Мандельбротом в частной беседе: фрактал как физическое тело – объект, в котором присутствуют дыры всех размеров. Действительно, средняя плотность головки швейцарского сыра меньше средней плотности отрезанного от неё кусочка: в нём и дыр меньше, и размеры их меньше.

Из равенства
$$\rho L^3 = \rho_c L^{D_m}$$
следует, что показатель $D_m$, равный

$$D_m = 3 - \ln(\rho_c/\rho)/\ln L,$$

определяет массовую фрактальную размерность через истинную и среднюю плотности сахара (что особенно удобно в экспериментальных исследованиях), и характерные размеры тела.

Массовую размерность $D_M$ можно определить и по-другому. Пусть нужно рассчитать массу фрактального шара. В этом случае зависимость массы от радиуса ведёт себя как

$$M(R) \sim R^{D_M},$$

где $D_M$ – массовая фрактальная размерность. Очевидно, что в пределе

$$D_M = \lim_{R \to \infty} (\ln M(R)/\ln R).$$

Так определённые массовые размерности являются глобальными характеристиками.

Для строго самоподобных математических фракталов, например, ковра Серпинского или канторовой пыли, массовая размерность $D_M$ совпадает с размерностью Хаусдорфа-Безиковича, потому что определяется размерностью подобия скейлингового закона, которая задаётся алгоритмом, порождающим фрактал [29].

**Химическая размерность**

Если рассматривать перколяционный кластер как решёточную модель ветвящейся полимерной молекулы, то узлы решётки при этом соответствуют мономерам, расстояния между узлами – химическим связям, а число шагов по кластеру – количество химических связей вдоль пути по молекуле от $i$-го мономера к $j$-му – можно определить как химическое расстояние. Размер-ность, связанная с этим расстоянием, называется химической, или размерностью связности $D_{ch}$ [6].

Для определения $D_{ch}$ рассмотрим шар $B_{ch}$ как множество узлов, для которых $R_{ch} \leq n$, и определим химическую размерность как показатель такой, что количество узлов $N$, принадлежащих $B_{ch}$, растёт как $N \sim n^{D_{ch}}$.

Величина $D_{ch}$ есть отношение двух размерностей – фрактальной размерности кластера $D$ и размерности $D_R$ кривой, длина которой определяет химическое расстояние. Для двумерного случая численно определённое значение $D_{ch}$ равно 1,72, следовательно, размерность «геодезической» равна $D_R = D / D_{ch} \approx 1{,}10$; видно, что это «не слишком изломанная линия» [6].

**Эффективная размерность**

Эффективная размерность – понятие, которое выражает соответствие между математическими множествами и модельными объектами, и которому, по мнению Б. Мандельброта, не следует давать точного определения [20].

Как составляющая модельного описания, эффективная размерность обладает «особым взглядом». Известно, что макроскопические объекты, даже такие «тщедушные» как крылышки пчелы, семена клубники, шёлковая нить и паутина являются трёхмерными телами. Но в математических моделях можно полагать, что крылышки имеют размерность два, что размерность семян – нуль, нити – один, а паутины – между 1 и 2.

Субъективная составляющая эффективной размерности хорошо видна в примере Б. Мандельброта [20]. «Пусть есть шар диаметром 10 см, скрученный из толстой нити диаметром 1 мм. Удалённому наблюдателю клубок покажется фигурой с нулевой размерностью, т.е. точкой. С расстояния в 10 см шар из нитей выглядит как трёхмерное тело, а с расстояния в 10 мм – как беспорядочное переплетение одномерных нитей. На расстоянии в 0,1 мм каждая нить превратится в толстый канат, а вся струк-тура целиком опять станет трёхмерным телом. На расстоянии 0,01 мм «канаты» превратятся в переплетение волокон, и шар снова станет одномерным, и так далее. Наконец, когда клубок превратится в скопление, состоящее из какого-то конечного числа точек, имеющих размеры, сравнимые с атомными, его размерность снова станет равной нулю» [20].

**Вместо заключения**

В обзоре описаны некоторые универсальные и специальные размерности, вошедшие в математический аппарат теоретической физики в двадцатом столетии, как один из инструментов исследования физических систем, который позволяет расширить понимание многих явлений и процессов, иногда – сформулировать новые понятия и модели, кроме того, нередко приводит к неожиданно удачному и точному описанию, демонстрируя *непостижимую эффективность* [30] использования размерностей и в математике, и в естественных науках.

**Литература**

1. Poincaré H. Pourquoi l'espace à trois dimensions. // Revue de métaphysique et de morale. – 1912. – V. 20. – P. 483-504; перевод: Пуанкаре А. Почему пространство имеет три измерения. / В кн. А. Пуанкаре О науке. – М.: Наука, 1990. – С. 555-579.
2. Lebesgue H. Sur la non applibilité de deux domaines appartenant à des espaces à *n* et *n + p* dimensions (extrait d'une lettre à M.O. Blumental). // Math. Ann. – 1911. – V. 70. – P. 166-168.
3. Brouwer L.E.J. Über den natürlichen Dimensionsbegriff. // Journal Für Die Reine Und Angewandte Mathematik. – 1913. – V. 142. – P. 146-152.
4. Математический энциклопедический словарь. – М.: СЭ, 1988. – 848 с.
5. Гуревич В., Волмэн Г. Теория размерности. – М.: ГИИЛ, 1948. – 232 с.
6. Эфрос А.Л. Физика и геометрия беспорядка. – М.: Мир, 1982. – 176 с.

7. Соколов И.М. Размерности и другие критические показатели в теории протекания. // УФН. – 1986. – Т. 150, вып. 2. – С. 221-255.
8. Меньшиков М.В., Молчанов С.А., Сидоренко А.Ф. Теория перколяции и некоторые приложения. / Итоги науки и техники. – Серия Теория вероятностей. Математическая статистика. – Т. 24. – М.: ВИНИТИ, 1986. – С. 53-110.
9. Козлов С.М. Геометрические аспекты усреднения. // УМН. – 1989. – Т. 44, вып. 2. – С. 79-120.
10. Меньшиков М.В. Оценки перколяционных порогов для решеток в $R^n$. // Доклады АН СССР. – 1985. – Т. 284, № 1. – С. 36-39.
11. Жиков В.В. Асимптотические задачи, связанные с уравнением теплопроводности в перфорированных областях. // Математический сборник. – 1990. – Т. 181, вып. 10. – С. 1283-1305.
12. Физический энциклопедический словарь. – М.: СЭ, 1984. – 944 с.
13. Схоутен Я.А. Тензорный анализ для физиков. – М.: Наука, 1965. – 456 с.
14. Гринченко В.Т., Мацыпура В.Т., Снарский А.А. Введение в нелинейную динамику. Хаос и фракталы. – М.: ЛКИ, 2010. – 280 с.
15. Пуанкаре А. Наука и гипотеза. / В кн. А. Пуанкаре О науке. – М.: Наука, 1990. – С. 38-78.
16. Ehrenfest P. In that way does it becomes manifest in the fundamental laws of physics that space has three dimensions? // Proceedings of Royal Netherlands Academy of Arts and Sciences. – 1918. – V. 20, iss. 1. – P. 200-209; reprinted in Paul Ehrenfest Collected Scientific Papers. / Ed. by M.J. Klein. – Amsterdam: North Holland Publ. Co., 1959. – P. 400-409.
17. Пуга В.А. Мультиразмерное гравитационное взаимодействие. Кривые вращения галактик. // ЖЭТФ. – 2014. – Т. 146, вып. 3 (9). – С. 500-512.
18. Турышев С.Г. Экспериментальные проверки общей теории относительности: недавние успехи и будущие направления исследований. // УФН. – 2009. – Т. 179, вып. 1. – С. 3-34.
19. Федер Е. Фракталы. – М.: Мир, 1991. – 254 с.

20. Мандельброт Б. Фрактальная геометрия природы. – М.: ИКИ, 2002. – 656 с.
21. Колмогоров А.Н., Тихомиров В.М. Е-энтропия и ε-ёмкость множеств в функциональных пространствах. // УМН. – 1959. – Т. 14, вып. 2 (86). – С. 3-86.
22. Гелбаум Б., Олмстед Дж. Контрпримеры в анализе. – М.: Мир, 1967. – 251 с.
23. Kahane J.P., Salem R. Ensembles parfaits et series trigonometriques. – Paris: Hermann, 1963. – 192 p.
24. Pontrjagin L., Schnirelman L. Sur une propriété métrique de la dimension. // Ann. Math. – 1932. – V. 33. – P. 156-162; перевод Понтрягин Л., Шнирельман Л. Об одном метрическом свойстве размерности. / В кн. [5], с. 210-218.
25. Божокин С.В., Паршин Д.А. Фракталы и мультифракталы. – Ижевск: НИЦ РХД, 2001. – 128 с.
26. Rényi A. Probability theory. – Amsterdam: North-Holland Publ. Co., 1970. – 666 p.
27. Кузнецов С.П. Динамический хаос. – М.: Физматлит, 2006. – 356 с.
28. Зельдович Я.Б., Соколов Д.Д. Фрактали, подобие, промежуточная асимптотика. // УФН. – 1985. – Т. 146, вып. 7. – С. 493 -506.
29. Шрёдер М. Фракталы, хаос, степенные законы. – Ижевск: НИЦ РХД, 2001. – 528 с.
30. Вигнер Е. Непостижимая эффективность математики в естественных науках. // УФН. – 1968. – Т. 94, вып. 3. – С. 535-546.

# СОДЕРЖАНИЕ